\documentclass[aap,noinfoline,reqno]{imsart}

\usepackage[total={5.8in, 8in},marginparwidth=2.5cm, marginparsep=.5cm]{geometry}
\usepackage{times}
\usepackage{soul}
\usepackage{currfile}
\usepackage[utf8]{inputenc}
\usepackage[english]{babel}
\usepackage{amsmath}
\usepackage{amsfonts}
\usepackage{amssymb}
\usepackage{amsthm}
\usepackage{mathtools}
\usepackage{graphicx,import}
\usepackage{comment}
\usepackage{marginnote}
\usepackage{xargs}                      
\usepackage[pdftex,dvipsnames]{xcolor}
\usepackage{microtype}
\usepackage{sub_JP}
\usepackage{mhequ}
\usepackage{imakeidx}
\usepackage{enumerate}
\makeindex

\startlocaldefs

\graphicspath{ {./} }
\numberwithin{equation}{section}
\DeclareMathAlphabet{\mathbbold}{U}{bbold}{m}{n}
\DeclarePairedDelimiter\abs{\lvert}{\rvert}%
\DeclareMathAlphabet\mathbfcal{OMS}{cmsy}{b}{n}

\newcommand{\figuresfolder}{.}

\newcommand{\h}[2]{h_r^{(#1)}(#2) }
\def\nablavU{\nabla_r^{(v)} U(x)}

\def\UQv1{L_r^{(v{{{}}})}(\xvi)}

\def\BB{\mathcal B}

\def\xvi{x{{{{}}}}}
\def\vvi{v{{{{}}}}}
\def\xvij{x{(j,k)}}
\def\vvij{v{(j,k)}}
\def\dsr{d_\star}
\newcommand{\Ov}{\pc{v^{-1}\mathbb N_0}^d}
\newcommand{\Ovp}{\pc{v^{-1}\mathbb N}^d}
\newcommand{\Ovm}{\pc{v^{-1}\mathbb N_0}^m}
\newcommand{\qrq}{Q_r^{(\vvi)}(\xvi)}
\newcommand{\qqq}{Q_{r'}^{(\vvi)}(\xvi)}

\newcommand{\prrr}{P^{(\vvi)}_{r,r'}(\xvi)}

\newcommand{\twf}{Theorem~\ref{t:wf3}}
\newcommand{\eref}[1]{(\ref{#1})}
\newcommand{\fref}[1]{Fig.~\ref{#1}}
\newcommand{\tref}[1]{Theorem~\ref{#1}}
\newcommand{\lref}[1]{Lemma~\ref{#1}}
\newcommand{\aref}[1]{Assumption~\ref{#1}}
\newcommand{\dref}[1]{Def.~\ref{#1}}
\newcommand{\rref}[1]{Remark~\ref{#1}}
\newcommand{\cref}[1]{Cor.~\ref{#1}}
\newcommand{\pref}[1]{Proposition~\ref{#1}}
\newcommand{\exref}[1]{Example~\ref{#1}}

\newcommandx{\unsure}[2][1=]{\todo[linecolor=red,backgroundcolor=red!25,bordercolor=red,#1]{#2}}
\newcommand{\pc}[1]{\left(#1\right)}
\newcommand{\pq}[1]{\left[#1\right]}
\newcommand{\pg}[1]{\left\{#1\right\}}
\newcommand{\pd}[1]{\abs*{ #1}}
\newcommand{\eg}{\emph{e.g.}~}

\newcommand{\la}{\lambda}
\newcommand{\R}{\mathcal R}

\newcommand{\Ex}[2]{\mathbb E_{#1}\left[ #2\right]}

\newcommand{\px}[2]{\mathbb P_{\kern -0.2em #1}\left[ #2\right]}
\newcommand{\p}[1]{\mathbb P \left[ #1 \right]}

\newcommand{\RR}{\mathbb R_+}
\newcommand{\Rr}{\mathbb R}
\newcommand{\I}{\mathbb I}

\newcommand{\C}{\mathcal C}

\newcommand{\ie}{{\it i.e.}, }
\def\S{\mathcal S}

\newcommand{\Lv}{\mathcal L_v}

\newcommand{\net}{$(\S,\C,\R)$~}

\newcommand{\tg}[1]{\theta{{{}}}^{\langle w{{{}}} , #1 \rangle}}

\newcommand{\supp}{{\rm supp}}
\newcommand{\super}{{\rm super}}
\newcommand{\inn}{{\rm in}}
\newcommand{\outt}{{\rm out}}
\newcommand{\PP}{\mathcal P}
\newcommand{\AN}{\mathcal A_{\rho,\rho'}}

\newcommand{\VV}{\mathcal V}
\newcommand{\onevec}{1\!\!1}
\newcommand{\Ind}[1]{\I_{#1}}
\newcommand{\vp}{v'}
\newcommand{\anddd}{\qquad {\rm and } \qquad}
\newtheorem{theorem}{Theorem}[section]
\newtheorem{remark}[theorem]{Remark}
\newtheorem{df}[theorem]{Definition}
\newtheorem{lm}[theorem]{Lemma}
\newtheorem{ex}[theorem]{Example}
\newtheorem{coro}[theorem]{Corollary}
\newtheorem{cla}[theorem]{Claim}
\newtheorem{axx}{Assumption}

\newtheorem{pro}[theorem]{Proposition}

\newcommand{\abbr}[1]{{\small\sc\lowercase{#1}}}

\newcommand{\dfn}[1]{\begin{df} #1 \end{df}}
\newcommand{\prp}[1]{\begin{pro} #1 \end{pro}}

\newcommand{\rmk}[1]{\begin{remark} #1 \end{remark}}
\newcommand{\lmm}[1]{\begin{lm} #1 \end{lm}}
\newcommand{\trm}[1]{\begin{theorem} #1 \end{theorem}}
\newcommand{\exm}[1]{\begin{ex} #1 \end{ex}}

\def\Label#1{}

\newenvironment{myenum}
{\begin{enumerate}[(a)]
  \setlength{\itemsep}{1pt}
  \setlength{\parskip}{0pt}
  \setlength{\parsep}{0pt}}
{\end{enumerate}}
\def\d{{\rm d}}
\def\dt{{\d t}}

\def\emx#1{\emph{#1}\index{#1}}
\let\epsilon=\varepsilon
\let\rho=\varrho
\let\kappa=\varkappa
\def\xt{\zeta(t)}
\def\x{\zeta}
\endlocaldefs

\makeatletter
\renewcommand*{\@fnsymbol}[1]{\ensuremath{\ifcase#1\or *\or \dagger\or
\mathsection\or
     \dagger \or \ddagger\or  \|\or **\or \dagger\dagger
     \or \ddagger\ddagger \else\@ctrerr\fi}}
\makeatother

\begin{document}

\begin{frontmatter}

\title{Large deviations theory for Markov jump models of chemical reaction networks}
\runtitle{Large deviations theory for chemical reaction networks}

 \author{\fnms{Andrea} \snm{Agazzi}\corref{AA}\ead[label=e1]{agazzian@gmail.com}\thanksref{t1,t2,t3}}
\thankstext{t1}{Research partially supported by NSF grant DMS-1613091.}
\thankstext{t2}{Research supported by ERC advanced grant 290843 and SwissMAP.}
\thankstext{t3}{Research partially supported by SNSF grant 161866 and SwissMAP.}

\author{\fnms{Amir} \snm{Dembo}\thanksref{t1}}
\and
\author{\fnms{Jean-Pierre} \snm{Eckmann}\thanksref{t2}}
\affiliation{Stanford University\thanksmark{t1} and University of Geneva\thanksmark{t2}}

\runauthor{A.~Agazzi, A.~Dembo, J.-P.~Eckmann}


\begin{abstract} \\
We prove a sample path Large Deviation Principle (\abbr{ldp})
for {a class} of jump processes whose rates
are not uniformly Lipschitz continuous in phase space. Building on it
we further establish the corresponding Wentzell-Freidlin (\abbr{W-F})
(infinite time horizon) asymptotic theory.
These results apply to jump Markov processes that model the dynamics
of chemical reaction networks under mass action kinetics, on a
microscopic scale.
We provide
natural sufficient topological conditions for the
applicability of our \abbr{ldp} and \abbr{W-F} results. This then justifies
the computation of non-equilibrium potential and exponential transition time
estimates between different attractors in the large volume
limit, for systems that are beyond the reach of standard
chemical reaction network theory.
\end{abstract}



\end{frontmatter}


\section{Introduction}

The dynamics of chemical reactions are usually modeled by mass-action
equations: A system of a polynomial ordinary differential equations which relate the evolution of concentrations of chemical compounds.
These systems of equations inherit their structure from the topology of the Chemical Reaction Network (\abbr{CRN}) they model, and the interplay between topology and dynamics of mass action systems is the object of study of chemical reaction network theory \cite{anderson11, feinberg87, horn72}.
These sets of \abbr{ODE}s approximate the interactions of the individual molecules involved. The discrete nature of chemical reaction systems can be captured by discrete models where the state of the system is given by the number of molecules of each type that are present in the reactor. In this framework, when a reaction occurs, the input molecules combine to form the output ones, and the system jumps to a new state. The dynamics of such systems are in general modeled stochastically as a pure jump Markov process \cite[Sec.~11, Example~C]{either86} whose jump rates are approximations of the reaction rates found in deterministic mass action models. Finally, assuming that the system has volume $v$, one can study how the stochastic dynamics of the process $X_t^v$ describing the concentration of the different chemical species at time $t$ scale with the parameter $v$. This is the object of study of this paper.

Similar discrete stochastic mass action kinetics models have been applied to disease propagation dynamics \cite{rock14}, genetic algorithms \cite{nobile16}, and for the simulation of noisy biochemical reaction networks through the application of the so-called Gillespie algorithm \cite{gillespie00}. Asymptotics such as limit theorems on the convergence of the stochastic trajectories towards the deterministic ones have been proven in the probability literature \cite{either86}. More recently, results on product-form steady state distributions for a certain class of \abbr{CRN}s have been obtained in {\cite{anderson15, cappelletti16}} and conditions for the irreducibility and ergodicity of the stochastic chemical dynamics of reaction networks have been presented in \cite{gupta13,pauleve14}. Our work extends these results to the domain of large deviations theory, identifying a large class of
\abbr{CRN}s to which that theory applies. We prove in
particular that Wentzell-Freidlin exit time estimates can be applied
to such systems, rigorously justifying the widespread use of potential
theory \cite{toth99,qian00,li1} and ultimately allowing for the analysis
of events that play a key role in, {\it e.g.}, theoretical
biochemistry \cite{balaszi11,qian10} {and that} are not covered by deterministic mass
action models, because deterministic models do not allow for
transitions between different attractors.

\subsection{The model and its sample path \abbr{ldp}}

We consider a set of \emx{chemical species} $\S = \{s_1, s_2, \dots  ,s_d\}$,
whose interactions are described by a finite set of \emph{reactions}\index{reaction} $\R = \{r_1, r_2, \dots , r_m\}$\,. {Throughout, we denote by $\mathbb N_0$ the set of natural numbers including $0$.}
Each reaction is uniquely identified by its \emph{substrates}\index{substrate} (input
species) and \emph{products}\index{product} (output species), and
we express such a reaction as $r = \{{ c}^r_{\rm in} \rightharpoonup c^{r}_{\rm
  out}\}$, with $c^{r}_{\rm out}$, ${ c}^r_{\rm in} \in\mathbb N_0^d $ representing the
multiplicity of the species {$s_i \in \S$} in the input or output of the reaction.
The set $\C$ of \emph{complexes}\index{complex} consists
of all $c_{\#}^r$ (with $\#=$ ``in'' or ``out''), and
for each reaction $r \in \R$ we
define the \emx{reaction vector} ${ c}^r := c^{r}_{\rm out}-{c}^r_{\rm
  in} \in \mathbb Z^d$.
A \abbr{CRN} is thus defined by the triple {\net.}

\exm{The system
\begin{equ}
\emptyset \mathrel{\mathop{\rightharpoonup}^{\mathrm{r_1}}} A+B \mathrel{\mathop{\rightharpoonup}^{\mathrm{r_2}}} 2B \mathrel{\mathop{\rightharpoonup}^{\mathrm{r_3}}} A
\label{e:ex1}
\end{equ}
is a \abbr{CRN} with $\S = \{A,B\}$ and $\R =
\{r_1,r_2,r_3\}$\,. The set of complexes of these reactions is $\C =
\{\emptyset,\{A+B\},\{2B\},\{A\}\} = \{(0,0),(1,1),(0,2),(1,0)\}$ (in the basis spanned by $(A,B)$).
\label{ex:1}}

In this paper, we study the behavior as a function of $v$ of the scaled process
\begin{equ}
(X_t^v)_i := v^{-1} (N_t)_i \,, \qquad i \in 1, \dots , d~,
\end{equ}
where $N_t \in \mathbb N_0^d$ represents the number of molecules of the $d$ species and $X_t^v \in \Ov$ denotes their number density (in mols) at time $t$.
The interactions among molecules are then described by each reaction
$r \in \R$ standing for a possible jump of the process
$ X_t^v \rightarrow X_t^v + v^{-1} { c}^r$, with ${c}^r$ the reaction (or jump) vector associated with $r \in \R$. Correspondingly,
$X_t^v$ is a continuous time pure jump Markov process with generator
\begin{equ}
(\mathcal L_v f) (x) := v \sum_{r \in \mathcal R}{\Lambda_r^{(v)}(x)} \left(f(x
{+} v^{-1} {c}^r)-f(x)\right)
\label{e:Lmarkov_1}
\end{equ}
for $f~:~\Ov \rightarrow \Rr$ and the volume-normalized {mass action kinetics} jump rates
\begin{equ}
\Lambda_r^{(v)}(x) =
k_r v^{-\|c_\inn^r\|_1} \prod_{i = 1}^d
\binom{v {x}_i}{ (c_\inn^r)_i} (c_\inn^r)_i!
\label{e:jrates}
\end{equ}
for some reaction (rate) constants $k_r>0$, where
$\binom{a}{b}$ denotes the binomial coefficient
which by convention is zero when $b \notin [0,a]$ and
$\|\cdot\|_1$ denotes the $\ell_1$-norm.
{The mean-field character of this model reflects the underlying assumption of homogeneous stirring of the reactor. The scaling in $v$ of the rate constants makes them asymptotically extensive quantities in \eref{e:Lmarkov_1} and takes into account that it is harder for molecules to meet as $v$ increases.}

\rmk{For a fixed volume $v$ and
initial condition $X_0^v = x_0^v \in \Ov$,
the process {$X_t^v$} is confined to
$S_{x_0^v}^v := \pg{x_0^v + \{\sum_{r \in \R} \alpha_r c^r~:~\alpha \in \Ovm  \}} \cap \RR^d$,
{where $\RR$ represents the set of nonnegative real numbers.}
Indeed, $X_t^v$ cannot jump outside of $\Ov$ since
$\Lambda^{(v)}_r(x)=0$ for any $r \in \R$ such that
$x + v^{-1} c^r \notin \Ov$ so the corresponding
summand in \eref{e:Lmarkov_1} is then zero
(regardless of $f(\cdot)$).
\label{r:zerosummands}
}



\rmk{{In the limit $v \rightarrow \infty$,
the sample paths of the processes $X_t^v$ starting at
$X_0^v = x_0^v \to x_0 \in \RR^d$ almost surely converge---{uniformly over $[0,T]$ for any $T > 0$\,}---to the solution $\xt $ of the deterministic \abbr{ODE}
\begin{equ}
\frac{\d \x}{\dt}  = \sum_{r \in \mathcal R} \lambda_r(\x) {c}^r \,, \qquad  \qquad \x(0) = x_0~,
\label{e:ma}
\end{equ}
having the asymptotic reaction rates
\begin{equ}
\lambda_r(x) := k_r \prod_{i = 1}^d {x}_i^{(c^r_{\rm in})_i}
\label{e:asrates}
\end{equ}
{provided that a solution of \eref{e:ma} exists up to time T}
(see \cite[§11, Thm.~2.1]{either86}, where such a
a {Functional Law of Large Numbers} (\abbr{flln}) is derived for {certain} \abbr{crn}s).}
\label{r:boundednessoftrajectories}}

We show in Section \ref{s:ldp} that under the following mild assumption on
the generator $\Lv$ of the scaled process, the solution $X_t^v$ to the corresponding martingale problem satisfies a sample path \abbr{ldp}
in the supremum norm, with an explicit rate function (see \tref{t:largedev}).
{While proving this \abbr{ldp} we also verify that in this setting
the \abbr{ode} \eref{e:ma} admits global solutions (and that
the \abbr{flln} of \rref{r:boundednessoftrajectories} holds).}
\begin{axx} Let $X_t^v$ be the solution of the martingale problem generated by the generator $\Lv$ of \eref{e:Lmarkov_1}\,. We assume
\begin{myenum}
\item There exist $b<\infty$ and a continuous, positive
function $U(x)$ of compact level sets, such that for
some non-decreasing function $\vp :\RR \to \RR$,
\begin{equ}
(\Lv U^v) \, (x) \leq
e^{b v}
\qquad \forall v > \vp (\|x\|_1),\;\; x \in (v^{-1} \mathbb N_0)^d ~,
\label{e:expdrift}
\end{equ}
{where $U^v(\cdot)$ denotes the $v$-th power of $U(\cdot)~$.}
\item
{With positive probability, starting at $X_0^v=0$ the Markov
process $X_t^v$ reaches in finite time some}
{state $x_+$ in the strictly positive orthant $(v^{-1} \mathbb N)^d$.}
\end{myenum}
\label{a:1}
\end{axx}

\rmk{The existence of a solution $X_t^v$ to the martingale problem generated by $\Lv$ with initial condition $x_0^v \in \Ov$ is guaranteed by standard theory (see \cite[Thm.~8.3]{bellet06}), up to the possibility of explosion.
In \lref{t:exptight} we show that this possibility is ruled out by \aref{a:1}.}

 \rmk{\aref{a:1}(b) requires that all chemical species can be created, at least indirectly, starting {from zero, hence from any other}
possible state of the system. In particular,
there must exist at least one chemical reaction without substrates,
namely, with $c^r_\inn=0$. Such {constant rate}
reactions
are used, \emph{e.g.}, in {mass action} models of cellular dynamics
\cite{anderson15} and continuous-flow stirred-tank chemical reactors
\cite{feinberg87}, {to}
model inflow of chemicals from the environment (correspondingly, these \abbr{crn}s often also have certain products exit the network, reflected by a mass
loss in some reactions).
It is possible to have an \abbr{LDP}
without \aref{a:1}(b), {but then}
 even when starting at $x_0^v \to x_0$ which is strictly positive, we
 may have a path of finite rate that leads to $\partial \RR^d$ and
 stays there forever. This would create problems establishing the
 Wentzell-Freidlin estimates.\label{r:ired-in}}
Proceeding to state our sample path \abbr{ldp}, hereafter
$D_{0,T}\pc{\RR^d}$ denotes the
space of \emph{c\`adl\`ag} functions $z~:~[0,T]\rightarrow \RR^d$
equipped with the topology of {uniform convergence}.
For $z(\cdot)$ in the subspace $AC_{0,T}\pc{\RR^d}$
of absolutely continuous functions from $[0,T]$ to $\RR^d$,
let $z'(\cdot)$ denote its Radon-Nikodym derivative
with respect to Lebesgue measure. Further,
for $\lambda = (\lambda_r) \in { \RR^m}$, {$q=(q_r) \in \RR^m$,}
$\xi \in \Rr^d$ and $c^r \in \Rr^d$, let

\begin{align}
L(\lambda,\xi) &:= \sup_{\theta \in \Rr^d} \Big\{ \langle \theta, \xi \rangle  -
\sum_{r \in \R} \lambda_r \big[\exp \pc{\langle \theta, {c}^r \rangle}-1\big]\Big\}
\nonumber \\
&=  \inf \Big\{ \sum_{r \in \R}
\big[\lambda_r - q_r+ q_r \log \frac{q_r}{\lambda_r} \big]
~:~{q} \in Q_\R(\xi)\Big\} \,,
\label{e:lagrangian}
\end{align}
where
$Q_\R(\xi) := \{{q \in \RR^m} :~\sum_{r \in \R} q_r c^r = \xi\}$ and
$\langle \theta, \xi \rangle$ is the inner product of
$\theta,\xi \in \Rr^d$.
\trm{
For $\lambda_r(\cdot)$ of \eref{e:asrates} and
any {$x^v_0 \to x_0 \in \RR^d$,} under \aref{a:1}
the
sample paths $\{X_t^v~:~t \in [0,T]\}$
{with $X_0^v=x_0^v$,} satisfy the \abbr{LDP}
in $D_{0,T}\pc{\RR^d}$ with rate $v$ and the good rate function
\begin{equs}
I_{x_0,T}(z) := \begin{cases}
\int_0^T L\pc{\lambda(z(t)),z'(t)} \d t \quad &\text{if}\ z(0) = x_0\ \& \ z \in AC_{0,T}\pc{\RR^d},\\
\infty \quad &\text{otherwise}~.
\end{cases}
\label{dfn:Ixt}
\end{equs}
That is, for any set $\Gamma \subset D_{0,T}(\RR^d)${, denoting by $\Gamma^o$ and $\bar \Gamma$ the interior and, respectively, the closure of $\Gamma$,} we have
\begin{align}
\label{dfn:ubd}
\limsup_{v \rightarrow \infty} \frac{1}{v} \log \px{x_0^v}{X_t^v \in \bar \Gamma} &\leq - \inf_{z \in \bar \Gamma}I_{x_0,T}(z) ~,\\
\liminf_{v \rightarrow \infty} \frac{1}{v} \log \px{x_0^v}{X_t^v \in \Gamma^o} &\geq - \inf_{z \in \Gamma^o}I_{x_0,T}(z) ~.
\label{dfn:lbd}
\end{align}
\label{t:largedev}
}
\rmk{The identity \eref{e:lagrangian} is well known
(see \cite[Thm.~5.26]{LDT2}),
{and since the function $[b-u+u\log(u/b)]$ is positive
whenever $u \ne b$, it yields that the Lagrangian
$L(\lambda,\xi)$}
{vanishes
iff $\xi=\sum_{r \in \R} \lambda_r c^r$.}
Thus, the rate
$I_{x_0,T}(z)$
of \eref{dfn:Ixt} is zero iff $z(\cdot)$ solves {on $[0,T]$,}
the \abbr{ode}
\eref{e:ma} starting at $z(0)=x_0$ {(see \cite[Exercise 5.14]{LDT2}).}
\label{r:zero-rate}
}

\subsection{Topological stability and strongly endotactic networks}

Standard large deviations theory is not directly applicable for
proving \tref{t:largedev}, because we need to deal with jump rates that
are neither bounded away from zero, nor globally Lipschitz continuous.
The diminishing jump rates at the boundary are handled by adapting
our system to the framework of mean-field interacting particle systems,
and thereby applying \cite[Thm.~3.9]{dupuis16}, whereas
\lref{t:exptight} takes care of the lack of global Lipschitz continuity
by employing Lyapunov stability theory to establish exponential tightness.
In doing so, a most important challenge is to phrase a stability
condition strong enough for such exponential tightness,
and a sufficient condition for escape from the boundary (in extension
of \cite{SW2}), that are both applicable to a broad collection of
\abbr{CRN}s.

This is precisely what we do next, with our topological conditions
summarized by \aref{a:ase} below. Specifically, given a finite set
$Q \subset \Rr^{d}$ and a vector $w \in \Rr^{d}$, we call
\begin{equs}\Label{dfn:w-max}
Q_w := \{c \in Q~:~ \langle w, c - c' \rangle \ge 0 \ \text{ for all } c'
\in Q\} ~\,,
\end{equs}
the $w$-maximal subset of $Q$ and consider the following collection of
\abbr{crn}s.
\dfn{{\rm \cite{gopal13}} The network $(\S,\C,\R)$ is called \emph{strongly
    endotactic}\index{strongly endotactic} if
for any {non-zero} $w\in\Rr^d$, the set $\R_{w} \subseteq \R$ of
reactions such that $c_\inn^r \in (\C_\inn)_w$ contains at least one
reaction satisfying $\langle w, c^r \rangle < 0$ and no reaction
with $\langle w, c^r \rangle >  0$.
\label{d:Pendo1}}

\noindent
This class of \abbr{CRN}s is well known
(see \cite{gopal13}), and algorithms to determine if a network is strongly endotactic are devised in \cite{johnston14}
(using variants of the simplex algorithm).

\begin{figure}
\centering
\def\svgwidth{.4\textwidth}
\includegraphics[scale=0.75]{\figuresfolder/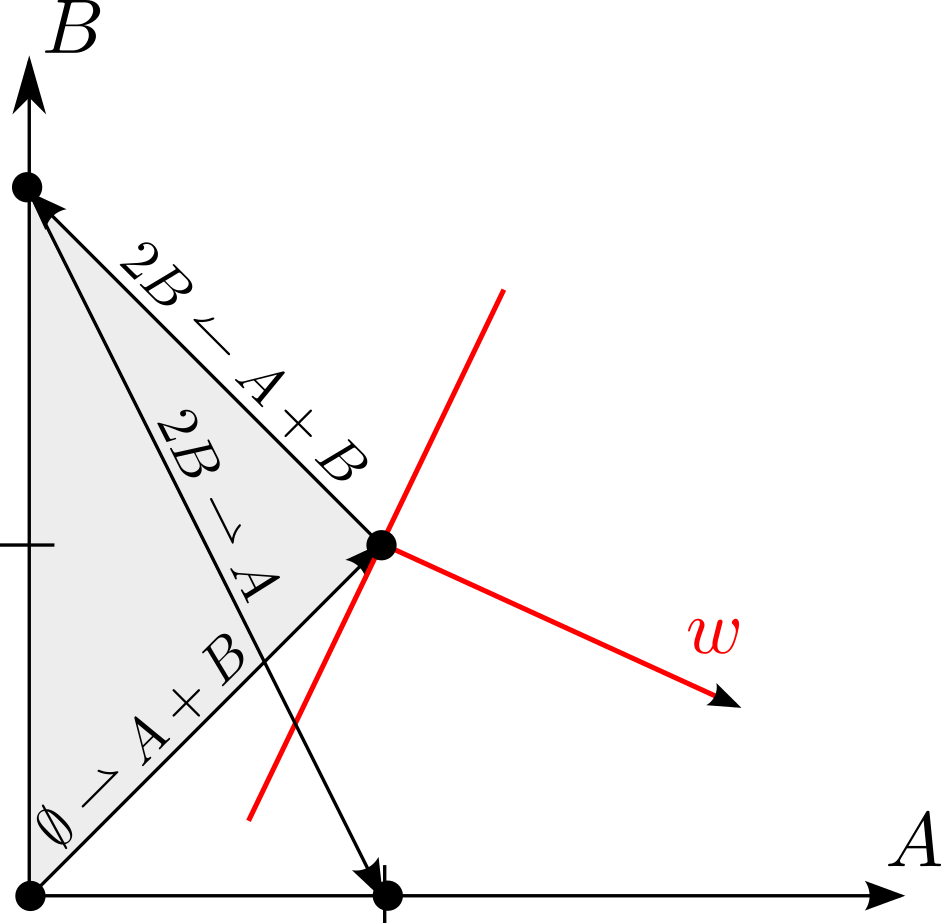}
        \caption{The $\RR^d$-diagram for \exref{e:ex1}. A vector $w$ in the space of complexes and the corresponding orthogonal hyperplane has been drawn in red to identify the $w$-maximal subset of the input complexes $(\C_\inn)_w$: the complex $A + B$\,.}
        \label{f:ex2}
\end{figure}

\addtocounter{theorem}{-8}
\noindent
\exm{[continued] {The network of \exref{e:ex1} is
represented as in \fref{f:ex2}, where we
identify $(\C_\inn)_w$ by sweeping $\RR^d$ with a hyper-plane orthogonal
to $w \in \Rr^d$ (here for $d=2$, drawn in red), and taking the last point of
$\C_\inn$ that such hyper-plane intersected. It is easy to see that our
specific network satisfies the requirements of \dref{d:Pendo1} and is therefore strongly endotactic.}}
\addtocounter{theorem}{7}

While in a strongly endotactic reaction network, all reactions ``point inward'' with respect to the faces of the convex hull of $\C_\inn$ (etymologically \emph{endo-tactic}: inward-arranged), our \abbr{LDP} requires addressing
the following additional boundary concept.
\dfn{A {non-empty}
subset {$\PP \subseteq \S$} is called a \emx{siphon} if every reaction
$r \in \R$ with {at least one output from $\PP$ also has
some input species from $\PP$.}}

\exm{It is readily checked that the sets $\PP = \{A\},\{A,B\}$ are siphons of the network
\begin{equ}
A \rightharpoonup 2A \rightharpoonup 3A+2B \rightharpoonup A~,
\Label{e:ex32}
\end{equ}
whereas $\{B\}$ is not.
}

We make the following assumption on the topological structure of \abbr{CRN}s.
{We call \net\ an
\emph{Asiphonic Strongly Endotactic} (\abbr{ASE}) network if it satisfies
\begin{axx}
The \abbr{CRN} \net has the properties:
\begin{myenum}
\item It is strongly endotactic, as in \dref{d:Pendo1},
\item It has no siphon $\PP \subseteq\S$\,.
\end{myenum}
\label{a:ase}
\end{axx}}

\rmk{\aref{a:ase}(b) is equivalent to {finding}, for any non-empty
$\PP \subseteq \S$, some reaction from $\R$ that produces at least one
output in $\PP$ while requiring no input species from $\PP$. When this
holds, then, for any state $x$ on the $\PP$-boundary of $\RR^d$ (namely,
having $x_i=0$ for all $s_i \in \PP$), there is some reaction of
non-vanishing rate
that brings the system back to a higher-dimensional subspace of $\RR^d$.
Following a sequence of such jumps we conclude that any
asiphonic \abbr{CRN} satisfies \aref{a:1}(b). {This definition coincides with the one of \emph{exhaustive} \abbr{crn}\emph{s} introduced in \cite{kammash14}.}\\
\label{r:irred}}

Combining the following result with \rref{r:irred} yields the \abbr{LDP} of
\tref{t:largedev} for the \abbr{ASE} networks of \aref{a:ase}.
\prp{[Existence of a Lyapunov function]
If the network is \abbr{ASE}, the generator $\Lv$ of \eref{e:Lmarkov_1}
satisfies \aref{a:1}(a)
{for the \emph{chemical Lyapunov}
(continuous) \emph{function }
\begin{equ}
U(x) :=  d + 1 + \sum_{i = 1}^d {x}_i (\log {x}_i - 1) \, : \RR^d \to \Rr_{\ge 1} \,.
\label{e:lyapunov1}
\end{equ}
\label{t:thm}}
}
{The connection between Lyapunov stability analysis and large deviations rate functions is an active area of research (see for example \cite{budhiraja15}).} Also, the problem of stability of mass action kinetics systems has been
addressed in \cite{anderson15, feinberg87, gopal13, horn72} and sufficient conditions for the existence of a globally attracting steady
state for the deterministic dynamics of such systems have been established in {\cite{anderson11,craciun15,feinberg87}}. In particular, the existence of a global attractor for a certain class of \abbr{CRN}s is proven
in {\cite{anderson11,gopal13}}
using the chemical Lyapunov function of \eref{e:lyapunov1}.
These results have been extended in \cite{gopal13} where the same
function is used for showing the existence of a compact attracting set for
strongly endotactic \abbr{crn}s. However, none of the references above deal directly with the generator $\Lv$, using the chemical Lyapunov function to establish exponential tail estimates for the finite-time distributions of such stochastic processes, as we do
in Section~\ref{s:prop} (where we prove Proposition \ref{t:thm} by verifying \eref{e:expdrift} for this function).

{\rmk{Proposition~\ref{t:thm} implies that it is sufficient to check a
    set of {graphical} conditions to guarantee the applicability of a
    \abbr{LDP} to the dynamics of \abbr{crn}s. This is most
    advantageous for applications in \emph{e.g.}, biochemistry, where
    typically $d>100$ and quantitative estimates like
    \eref{e:expdrift} would be prohibitive to check. {Note furthermore
    that our conditions do \emph{not} depend on the reaction rate
    constants $k_r$, which are often very difficult to determine.}}}

\subsection{Quasi-potential and exit time asymptotic}

Following the Wentzell-Freidlin approach, we utilize our sample path \abbr{ldp}
to define the corresponding
quasi-potential (as in \cite{fw98}), and provide
asymptotic analysis over an infinite time horizon, for
quantities of interest such as the exit time from some domain $\mathcal D \subset \RR^d$,
or the transition time between different attractors of \eref{e:ma}
(as proposed by \cite{toth99}). To do so, we first assume that
the domain of interest $\mathcal D$ has the following mild regularity properties.
\begin{axx}
The compact $\mathcal D \subset \RR^d$ is the closure of its interior, with boundary
$\partial \mathcal D$ that is a piecewise twice continuously differentiable
sub-manifold of $\RR^d$. Furthermore, there exists a
ball $\BB \subset \mathcal D$ so that for all $x \in \BB$ and $y \in \mathcal D$ the
set $\mathcal D$ contains the line segment between $x$ and $y$.
%
%
\label{a:2}
\end{axx}

\dfn{The \emph{quasi-potential} between any $x, y \in \RR^d$ is
\begin{equs}\Label{dfn:quasi}
\VV_{\mathcal D}(x,y) := \inf_{t \geq 0} \; \inf_{
\{z(\cdot) \in \mathcal D, z(t) = y\}} \; \{ I_{x,t}(z) \} ~\,,
\end{equs}
for $I_{{x},t}(z)$ of \tref{t:largedev}. We say that $x,y$
are $\mathcal D$-equivalent (denoted $x \sim_{\mathcal D} y$),
if $\VV_{\mathcal D}(x,y) = \VV_{\mathcal D}(y,x) = 0$. We further define
\begin{equ}
\VV_{\mathcal D}(A,B) := \inf_{x \in A, y \in B} {\VV}_{\mathcal D}(x,y)\,, \qquad
\forall A,B \subseteq \mathcal D \,,
\end{equ}
and use $\VV(\cdot,\cdot)$ for $\VV_{\RR^d}(\cdot,\cdot)$.
\Label{d:cost}}
{The equivalence $x \sim_{\mathcal D} y$ defines compact sets $K_i \subset \mathcal D$ where the
process can move {with probability
$\exp(-o(v))$.}
Throughout, we make the following assumption about their structure.}
\begin{axx}\emph{\cite[Condition A, §6.2]{fw98}}
{There exist $\ell$ compact sets $K_i \subset \mathcal D$ such that:
\begin{myenum}
\item every $\omega$-limit set of \eref{e:ma} lying entirely in $\mathcal D$ is
fully contained within one $K_i$\,,
\item for any $x \in K_i$ we have $x \sim_{\mathcal D} y$ if and only if $y \in K_i$\,,
\item {for all $K_i$ , the set $K_j$ minimizing $\VV_{\mathcal D}(K_i,K_j)$ is unique.}
\end{myenum}
We further assume that the conic hull ${\rm Co}\{{c}^r\}_{r \in \R}$
of vectors $\{{c}^r\}_{r \in \R}$ is $\Rr^d$.
\label{a:3}}
\end{axx}

%

Such $K_i$ are called \emph{stable} if $\VV(K_i^\delta,(K_i^\delta)^c)>0$
for $\delta>0$ small enough (where $B^\delta$ denotes the $\delta$-neighborhood
of the set $B$ in the $\|\cdot\|_1$-norm).
The most probable transitions between $\{K_i^\delta\}$ for small $\delta>0$
and $v \to \infty$ {(\ie those transitions that connect any such $K_i^\delta$
to the unique $K_j^\delta$ with $i \neq j$ minimizing $\VV_{\mathcal D}(K_i^\delta, K_j^\delta)$)}
define a deterministic dynamic on the finite
collection of stable compact sets. Such dynamic can be partitioned into
disjoint \emph{cycles}, with each cycle $\pi$ consisting of a single
transitive
point ($\pi=\{i\}$) or a periodic orbit $\pi =\{i_1\to i_2\to \dots\to
i_j\to i_1\}$ (c.f. \cite[§6.6]{fw98} for the precise definition).
Thanks to \aref{a:2} and \ref{a:3}, adapting the machinery of \cite{fw98}
to our setup, we transfer in Section \ref{s:wf} the sample path
\abbr{ldp} to the following result about the time it takes the
\abbr{crn} to exit $\mathcal D$ or a cycle $\pi$ and the probability
cost of relevant exit paths.

\trm{\emph{\cite[Thm.~5.1, 5.3 and 6.2, §6]{fw98}}

\noindent
Consider a \abbr{crn} satisfying \aref{a:1} and
the process $t \mapsto X_t^v$ starting at $x_0^v \to x \in \mathcal D^o$.
Let $\tau_{\mathcal D}$ denote its exit time from a set $\mathcal D$ that
satisfies \aref{a:2} and \ref{a:3}
and let
$\tau_\pi$ its first hitting time of $\cup_{j \notin \pi} K_j^\delta$
for a cycle $\pi \subseteq \{1,\ldots,\ell\}$ and sufficiently small
$\delta > 0$\,.
Then, with non-random $M_{\mathcal D}(x)$, $W_{\mathcal D}$ and $W_{\mathcal D}(x,y)$ as in \cite[§6]{fw98},
we have that for any $x$ in a compact $F \subset \mathcal D^o$ and $y \in \partial \mathcal D$
\begin{equ}
\lim_{\delta \to 0}
\lim_{v \rightarrow \infty}\frac{1}{v}\log\px{x_0^v}{
\| X_{\tau_{\mathcal D}}^v - y \|_1 < \delta} = W_{\mathcal D} -W_{\mathcal D}(x,y)\,,
\label{e:twf1}
\end{equ}
\begin{equ}
\lim_{v \rightarrow \infty} \frac{1}{v} \log \Ex{x_0^v}{\tau_{\mathcal D}} = W_{\mathcal D} - M_{\mathcal D}(x) \,.
\label{e:twf2}
\end{equ}
Furthermore, with
$C(\pi)<\infty$ as in \cite[§ 6.6]{fw98}, any $\gamma>0$ and
uniformly in $x \in \cup_{i \in \pi} (K_i)^{\delta/2}$,
\begin{equ}
\lim_{v \rightarrow \infty} \px{x_0^v}{\pd{ v^{-1}\log \tau_\pi - C(\pi)}\leq \gamma} = 1\,.
\label{e:twf4}
\end{equ}
\label{t:wf3}}

\rmk{Note that models in cell biology \cite{qian10} {usually have significantly larger
  dimension $d$
than many other applications of
Wentzell-Freidlin theory.}}


\section{Proof of \tref{t:largedev}}
\label{s:ldp}
We start by showing that \aref{a:1}(a) yields exponentially negligible
exit probability from the compact level sets of the function $U(\cdot)$.
\lmm{Let $\{X_t^v\}$ be a Markov jump process with generator \eref{e:Lmarkov_1} and initial condition $x_0^v \in \Ov$. Under
\aref{a:1}(a), there is, for every $\alpha$, $\beta$, $\gamma$, a finite
$\rho_{\alpha,\beta,\gamma}$, so that
\begin{equ}
\limsup_{v \rightarrow \infty}
\frac{1}{v}\log\Big(
{
\sup_{\|x_0^v\|_1 \le \gamma}} \,
\mathbb P_{x_0^v}\big[\sup_{t \in [0,e^{\beta v}]} {\|X_t^v\|_1} >
\rho_{\alpha,\beta,\gamma}
 \big]\Big) \leq - \alpha~.
\label{e:exptight}
\end{equ}
\label{t:exptight}}

\begin{proof} For each $\ell$ there is a $\rho=\rho(\ell)$ so that $\{x : U(x) \le \ell \}$
is a subset of the ball
\begin{equs}
\widetilde{K}_\rho := \{ x \in \RR^d : \|x\|_1 \le \rho \} \,.
\label{dfn:K-tilde}
\end{equs}
Considering the $v$-dependent stopping times
\begin{equ}\label{e:tau-r}
\sigma_\rho := \inf\{t > 0 ~:~ X_t^v \notin \widetilde{K}_\rho \}\,,
\end{equ}
and the stopped processes $\widehat X_t^{v,\rho} := X_{\sigma_\rho \wedge t}^v\,$, it follows by Markov's inequality that for any $T$,
\begin{align*}
\mathbb P_{x_0^v}\Big[\sup_{t \in [0,T]} \|X_t^v\|_1 > \rho \Big]
&=
\mathbb P_{x_0^v}\Big[\|\widehat X_{T}^{v,\rho}\|_1 > \rho \Big] \\
&\le
\mathbb P_{x_0^v}\left[
U\big(\widehat X_{T}^{v,\rho}\big) > \ell \right] \leq \ell^{-v} \mathbb E_{x_0^v}\left[
U^v\big( \widehat X_{T}^{v,\rho}\big)\right] \,,
\end{align*}
from which we get \eref{e:exptight} once we show that
\begin{equ}
{\sup_{\rho \ge \gamma}}
\limsup_{v \to \infty} \frac{1}{v} \log
{\sup_{\|x_0^v\|_1 \le \gamma, T \le e^{\beta v}} \,}
\mathbb E_{x_0^v}\left[U^v\big(\widehat X_{T}^{v,\rho}\big)\right] < \infty \,.
\label{e:markov}
\end{equ}
To this end, as $U(\cdot)$ is continuous,
${\sup_{\|x\|_1 \le \gamma}} \{U(x)\} \le e^{\kappa}$ for some
$\kappa=\kappa(\gamma) <\infty$. Further,
{when $\| X_0^{v} \|_1 \le \rho$,}
the Markov process $\widehat X_t^{v,\rho}$ has the generator
$\mathcal L_v$ of \eref{e:Lmarkov_1} restricted to $\widetilde{K}_\rho$ and is confined for any $v \ge 1$
to a compact $(\widetilde{K}_\rho)^{{\bar c}}$ with ${{\bar c}} := \sup_r \|c^r\|_1<\infty$.
Thus, combining Dynkin's formula \cite[§5.1]{dynkin65} with \aref{a:1}(a)
we find that for some $v_\rho \in [1,\infty)$, all $v>v_\rho$ and
{$\|x^v_0\|_1 \le \gamma \le \rho$,}
\begin{equ}
\mathbb E_{x_0^v}\left[U^v(\widehat X_{T}^{v,\rho})\right] \le U^v(x_0^{v}) +
\mathbb E_{x_0^v}\Bigg[\int_0^{\sigma_\rho \wedge T}
\mathcal (\Lv U^v)(X_s^{v})~\d s\Bigg]
  \leq e^{\kappa v} + T e^{bv} \,.
\label{e:dynkin}
\end{equ}
Considering for $T \le e^{\beta v}$ the limit as $v \to \infty$ of $v^{-1}$ times
the logarithm of \eref{e:dynkin} leads to \eref{e:markov}
and thereby concludes the proof.
\end{proof}

\rmk{\lref{t:exptight} can be alternatively proved by defining a super-martingale from the condition \eref{e:expdrift} on our generator, and applying \cite[Thm.~4.20]{LDT3} to it.}

{The Markov jump process $X^v_{tT}$ corresponds to the generator
of \eref{e:Lmarkov_1}, now with reaction constants $T k_r$
{for which \aref{a:1} continues to hold.} This changes
$\lambda(\cdot)$ of \eref{e:asrates}
to $T \lambda(\cdot)$, hence transforms $I_{x_0,T}(z(t))$
into $I_{x_0,1}(z(tT))$
(since $L(T \lambda,y) = T L(\lambda, T^{-1}y)$).
Thus, \abbr{wlog}, we take hereafter $T=1$ and} proceed to establish the exponential tightness of an
exponentially equivalent  process $\widetilde X_t^v$.
\lmm{
{Under \aref{a:1}(a),}
the $C_{0,1}(\RR^d)$-valued processes $\widetilde X_t^v$ obtained by
linearly interpolating the jump points of $t \mapsto X_t^v$, form an
exponentially tight family in the uniform topology, which
{for uniformly bounded $\|x_0^v\|_1$} is further
exponentially equivalent to $\{X_t^v\}$ in the uniform topology
on $D_{0,1}(\RR^d)$.
\label{l:new}}
\begin{proof}
For any consecutive jumps of $X_t^v$ occurring at (random) times $t_1<t_2$
we set
\[
\widetilde X_t^v := X_{t_1}^v + \frac{t-t_1}{t_2-t_1} (X_{t_2}^v-X_{t_1}^v) \,.
\]
Hence, $\|X_t^v - \widetilde{X}_t^v\|_1 \le v^{-1} {\bar c}$ for
finite ${\bar c} := \sup_r \|c^r\|_1$, all $t \ge 0$, and $v$, yielding
the exponential equivalence of
$\{\widetilde X_t^v \}$ and $\{X_t^v\}$ (in the uniform topology).
As for the exponential tightness of $\{\widetilde X_t^v\}$ in
$C_{0,1}(\RR^d)$, note that for any $t > s$,
\[
\|\widetilde X_t^v - \widetilde X_s^v \|_1 \leq v^{-1} {\bar c} N_{[s,t]} (X^v)  \,,
\]
where $N_{[s,t]}(X^v)$ counts the number of jumps by $X_\cdot^v$ in the time
interval $[s,t]$. Further, as $\Lambda_r^{(v)}(x) \le \lambda_r(x)$ for all
$x \in \RR^d$, we have for $\sigma_\rho$ of \eref{e:tau-r} and any $v \ge 1$
the monotone coupling $N_{[s,t]}(X^v) \le M_{[s,t]}^{\rho}$ on $[0,\sigma_\rho]$,
where $M^\rho$ is a Poisson process of intensity {$v \bar \Lambda^{\rho}$ and}
\[
\bar \Lambda^{\rho} := \sup_{x \in (\widetilde{K}_\rho)^{{\bar c}}} \,
\Big\{ \sum_{r \in \R} \lambda_r (x) \Big\} \,.
\]
In view of the Arzel\`a-Ascoli theorem and \lref{t:exptight}, it thus
suffices for the stated exponential tightness of $\{ \widetilde X_t^v \}$
to show that
\begin{equ}\label{e:mod-cont}
\lim_{\delta \to 0}
\limsup_{v \to \infty}
v^{-1} \log \p{\sup_{0 \le s \le t \le (s + \delta) \wedge 1} \,
\{ M_{[s,t]}^{\rho} \} \ge v \epsilon} = -\infty \,,\;\; \forall \rho < \infty,
\epsilon > 0 \,.
\end{equ}
To this end, by tail estimates for the Poisson($2\delta {v} \bar \Lambda^\rho$) law,
for any $\epsilon>0$ and $\rho < \infty$,
\[
\lim_{\delta \to 0}
\limsup_{v \to \infty}
v^{-1} \log \p{  \,  M^\rho_{[0,2\delta]} \ge v \epsilon} = -\infty \,.
\]
Further, if $|t-s| \le \delta$ and $n=[1/\delta ]$, then
\[
M^\rho_{[s,t]} \le {\max_{i=0,\dots,n-1}} \{M^\rho_{[i\delta,(i+2)\delta]} \}
=: \bar M^\rho_\delta
\,.
\]
Hence, applying the union bound for the maximum
$\bar M^\rho_\delta$
of $n$ identically distributed
Poisson($2\delta {v} \bar \Lambda^\rho$)
variables yields
\eref{e:mod-cont}, and thereby concludes the proof.
\end{proof}

Let $\mathcal M_1(\S_\star)$ denote the probability simplex over
{$\S_\star = \{\star\}\cup \S$} and
$c^r_\star :=  \langle \onevec, c^r \rangle =
\|c^r_{\rm out}\|_1 - \|c^r_\inn\|_1$
the number of molecules gained (or lost, if negative)
in each chemical reaction. For $\rho>0$ and
$\{ \lambda_r(\cdot) \}$ of \eref{e:asrates} {such
that \eref{e:ma} admits a solution} {$\zeta: [0,1] \mapsto \RR^d$
(\ie no blowup on $[0,1]$),}
we consider $\mu(t)$ satisfying the \abbr{ODE}
\begin{equs}\label{e:lln-dupuis}
\frac{\d \mu}{\dt}  = \rho^{-1} \sum_{r \in \mathcal R} \lambda_r(\rho \mu\vert_\S)
(-c_\star^r,{c}^r) \,, \qquad \mu(0) \in \mathcal M_1(\S_\star) \,,
\end{equs}
establishing a strictly positive lower bound on $\{ \mu(t) \vert_\S\}$ that
holds uniformly over {$ \|\mu(0)\vert_\S\|_1 \le \gamma/\rho < 1$} with arbitrary, fixed
$\gamma$ and all $\rho$ large enough. {This quantity is a rescaled projection on $\mathcal M_1(\S_\star)$
of the \abbr{ode} \eref{e:ma} with
initial condition $\|\zeta(0)\|_1 \leq \gamma$
{provided $\sup_t \|\zeta(t)\|_1 \leq \rho$.}
 In other words, adding a ``vacuum'' species $\{\star\}$\,, we map
$\zeta(t)$ onto $\mu(t)$, describing the dynamics a system conserving the total number of molecules. Note that $\mu(t)$ can be seen as the empirical measure of an \abbr{ips} in the limit of infinite number of particles.}

\lmm{\label{l:p314dupuis}
{Let \aref{a:1} hold and assume that \eref{e:ma} has a solution {for $t \in [0,1]$.}
Then, for any $\gamma > 0$ and for some $\rho_0(\gamma)$,} if
$\rho \ge \rho_0(\gamma)$ and $\rho \mu(0)\vert_\S \le \gamma$, the solution $\mu(t)$ of \eref{e:lln-dupuis} satisfies
$\mu(t) \in \mathcal M_1(\S_\star)$ for $t \in [0,1]$.\\
{Further, there exist $D \in \mathbb N$
and $b=b(\rho)>0$ such that for any such $\mu(t)$ we have}
\begin{equs}\label{e:bd-413}
{\mu_{s_i}(t) \geq b t^{D} \,, \qquad \qquad \forall t \in [0,1], \;\; i=1,\ldots,d \,.}
\end{equs}}
\begin{proof} Starting at $\langle \onevec, \mu(0) \rangle =1$,  it follows from the definition of $c^r_\star$ that
$\langle \onevec, \mu(t) \rangle = 1$ for all $t \ge 0$, with
the bijection
\begin{equs}\label{e:bijec-ode}
\xt  = \rho \mu(t) \vert_{\S} =: \Psi(\mu(t)) \,,  \qquad
\mu_\star (t) {:=} 1 - \rho^{-1} \|\xt \|_1 \,,
\end{equs}
between $\mu(\cdot)$ of \eref{e:lln-dupuis} and
{the assumed finite solution}
$\x(\cdot)$ of \eref{e:ma} {with $\|\x(0)\|_1 \leq \gamma$.} In particular,
$\x(0)=\Psi(\mu(0)) \in \widetilde{K}_\rho$ of \eref{dfn:K-tilde}
yields $\x(\cdot) \in \RR^d$ and
the condition $\rho \mu(0)\vert_\S \le \gamma$ translates
into $\|\x(0)\|_1 \le \gamma$. Our claim that $\mu(t)
\in \mathcal M_1(\S_\star)$ for $t \in [0,1]$ is thus just
\begin{equs}\Label{dfn:rho-0}
\rho_0(\gamma) :=
\sup_{\|\x(0)\|_1 \le \gamma} \, \sup_{t \in [0,1]} \|\xt \|_1 < \infty \,,
\end{equs}
which holds for $\rho_0(\gamma) \le 1 + \rho_{1,0,\gamma}$ of
\eref{e:exptight} (indeed, simply contrast the
\abbr{flln} of \rref{r:boundednessoftrajectories}
with the exponential decay in $v$ of probabilities from \lref{t:exptight}).

{Next, for any $\rho >0$ we multiply each reaction constant $k_r$ by
$\rho^{\|c^r_\inn\|_1-1}$ and \abbr{wlog} set hereafter $\rho=1$.
Identifying $s_j=j$,}
split the \abbr{rhs} of \eref{e:lln-dupuis}
at coordinate $i$ to a sum over reactions in $\R^i_+ := \{ r \in \R : c^r_i>0 \}$
and over those in $\R^i_- := \{ r \in \R : c^r_i<0 \}$. The contribution from
$\R^i_+$ is a polynomial $P_i(\cdot)$ in $\{\mu_1,\ldots,\mu_d\}$ of positive
coefficients (namely $k_r c^r_i$, $r \in \R^i_+$). Further, $c^r_i < 0$
requires $(c^r_\inn)_i \ge 1$ so the contribution of $\R^i_-$ is of the
form $\mu_i Q_i(\mu)$ for another polynomial $Q_i(\cdot)$ with positive coefficients.
Let $e(t):=\mu(t) \vert_\S - y(t)$, for the solution $y(t)$ of the modified \abbr{ODE}-s
\begin{equ}
\frac{\d y_i}{\dt} = P_i(y(t)) - C y_i(t) \,, \qquad  i=1,\ldots,d,
\qquad y(0)=\mu(0)\vert_\S \,,
\label{e:lowerbound}
\end{equ}
where
\[
C := {1+}
\max_{i \le d} \, \sup \{ Q_i(\mu) : \mu \in \mathcal M_1(\S_\star) \}  < \infty \,.
\]
Each $P_i(\cdot)$ is increasing
\abbr{wrt} the natural partial order on $\RR^d$, hence
\begin{equ}\Label{e:suchabound}
\frac{\d e_i}{\d t} + C e_i = P_i(y+e) - P_i(y) + \mu_i (C - Q_i(\mu)) \ge 0 ~,
\end{equ}
as long as $e(t)$ and $y(t)$ are both in $\RR^d$, with a strict inequality
as soon as {$\mu_i(t)>0$}. Hence, starting at $e(0)=0$ and
$y(0) \in \RR^d$, we establish \eref{e:bd-413} by showing that the same inequality holds if one substitutes the solution $y(\cdot)$ of \eref{e:lowerbound} to $\mu(\cdot)$, uniformly over
all $y(0) \in \RR^d$.
We achieve this goal by
utilizing \aref{a:1}(b) in at most $d$ steps,
to get that for some $D_k \in \mathbb N$ and $b_k > 0$,
\begin{equs}\label{e:k-step}
y_i(t) \ge b_k t^{D_k} \,, \qquad \forall t \in [0,1], \;
y(0) \in \RR^d ,\; i \in \S_k \uparrow \{1,\ldots,d\} \,.
\end{equs}
Specifically, starting at $\S_0=\emptyset$
let $\S_k = \S_{k-1} \cup \partial \S_k$ for
$$
\partial \S_{k} := \{j  \notin \S_{k-1} ~:~ \exists r \in \R,
\,(c_\outt^r)_j > 0\,\mbox{ and } \,\forall l \not \in \S_{k-1}\,,
(c_\inn^r)_l = 0\}~.
$$
In particular, $\partial \S_1$ consists of {all product species in reactions
with $c^r_\inn=0$} and from \aref{a:1}(b) we know that $\partial \S_1$
is non-empty (see \rref{r:ired-in}). Such a reaction with
$c^r_\inn=0$ and an output $i \in \partial \S_1$ contributes
to $P_i(\cdot)$ a positive constant term $k_{r,i} := k_r c^r_i$.
For $y \in \RR^d$ any other reaction may only increase $P_i(y)$, hence
\[
\kappa_1 := \inf_{i \in \partial \S_1}
\inf_{y \in \RR^d}\, \{ P_i(y) \} \, > 0 \,.
\]
Bounding the solution of \eref{e:lowerbound} from below
taking $\kappa_1$ instead of $P_i(y(t))$, and considering
the worst case $y_i(0)=0$, we deduce that for
$k=1$, $D_1=1$ and any $i \in \partial \S_k$,
\begin{equ}\label{eq:bdk}
y_i(t) \ge \kappa_{k} \int_0^t e^{-C (t-s)} s^{D_k-1} \d s
{\geq} b_k t^{D_k} \,,
\end{equ}
for some $b_k=\kappa_k g(C,D_k)>0$ {and $t \in [0,1]$.}
Increasing to $k=2$, observe that
if $\S_{k-1} \ne \S$ then by \aref{a:1}(b) there must be a reaction
$r$ that has at least one product
not from $\S_{k-1}$ while all of its substrates
are from $\S_{k-1}$. In that case, the non-empty set $\partial \S_k$
consists of {the} products of such reactions {that are not in $\S_{k-1}$} and
for any $i \in \partial \S_k$ a reaction $r = r_i \in \R$
of this type contributes to $P_i(y(t))$ a positive term of the form
$$
k_{r,i} \prod_{l \in \S_{k-1}} y_l(t)^{(c^r_\inn)_l}
\ge k_{r,i} (b_{k-1} t^{D_{k-1}})^{\ell_i} \,,
$$
for $\ell_i := \| c^{r_i}_\inn \|_1$, where we relied on already
having the bound
\eref{e:k-step} for $l \in \S_{k-1}$.
Setting
\[
D_k := 1 + D_{k-1} \max_{i \in \partial \S_k} \{\ell_i\} \,,\quad
\kappa_k := \min_{i \in \partial \S_k} \{ k_{r,i} b_{k-1}^{\ell_i} \} \,,
\]
recall that other reactions may only increase $P_i(y(t))$, hence
for $i \in \partial \S_k$ and $t \in [0,1]$,
\[
P_i(y(t)) \ge \kappa_k t^{D_k-1} \,.
\]
Exactly as we have done for $k=1$ and $D_1=1$, inserting such a lower bound
into \eref{e:lowerbound} and considering the worst case solution
($y_i(0)=0$), results with \eref{eq:bdk}. Further lowering $b_k$
to have the same bound extend also to all $i \in \S_{k-1}$ and proceeding if
necessary to $k=3$ and beyond
exhausts finally all of $\S$ after at most $d$ steps.
\end{proof}

\begin{proof}[Proof of \tref{t:largedev}]
Recall the Skorokhod $J_1$-topology on $D_{0,1}(\RR^d)$ which
is metrizable
by the coarsening of the sup-norm
\begin{equ}
d_{J_1}(z_1,z_2) := \inf_{\tau} \big\{
\|\tau\|_\star + \sup_{s \in [0,1]}  \|z_1(s)-z_2(\tau(s))\|_1 \big\} \,,
\label{e:dj}
\end{equ}
where $\|\tau\|_\star := \sup_{s \ne t}
\log \big\{|\tau(s)-\tau(t)|/|s-t| \big\}$ for strictly increasing
$s \mapsto \tau(s)$ with $\tau(0) = 0$, $\tau(1) = 1$. By
\lref{l:new} and the inverse contraction principle
of \cite[Corollary 4.2.6]{LDT1}, it suffices to
establish the \emph{weak} \abbr{ldp} for {$\{\widetilde X_t^v\}$}
in the metric space $(D_{0,1}(\RR^d),d_{J_1})$
(in this standard reduction we also rely upon \cite[Lemma 1.2.18]{LDT1}
to upgrade from weak \abbr{LDP} to full \abbr{ldp} before employing the
inverse contraction, and on \cite[Thm.~4.2.13]{LDT1} to transfer
the \abbr{ldp} in the uniform topology from
{$\{\widetilde X_t^v\}$} to {$\{X_t^v\}$).}
Next, consider the Markov jump process {$X_t^{v,\rho}$}
of generator \eref{e:Lmarkov_1} and volume-normalized jump rates
\begin{equs}\label{e:dfn-Lambda-rho}
\Lambda_r^{v,\rho}(x) := \Lambda_r^{(v)} (x) \,\Ind~\!
(\|x\|_1 \le \rho-v^{-1} c^r_\star)\,,
\end{equs}
{where $\Ind~\!(A)$ is the indicator function over a set $A~$.}
Taking ${{\bar c_\bullet}} := \sup_r \|c^r_\inn\|_1 \vee \|c^r_\outt\|_1$ and
\begin{equs}\label{dfn:rho-val}
\sup_{v \ge 1} \{ \|x_0^v\|_1 \} + {{\bar c_\bullet}} \le \rho \,,
\end{equs}
assures that $\{{X_t^{v,\rho}}, v \ge 1 \}$ is confined to $\widetilde{K}_{\rho}$
of \eref{dfn:K-tilde} and in view of \lref{t:exptight},
\begin{equs}\label{e:exp-good-apx}
\lim_{\rho \to \infty} \limsup_{v \to \infty}
v^{-1} \log \px{x_0^v}{{X^{v,\rho}_t} \not \equiv {X^{v}_t} } = - \infty \,,
\end{equs}
{where $X_t^{v,\rho} \not \equiv X_t^v$ represents the event where the paths of {$X_t^v$} and of {$X_t^{v,\rho}$} (coupled to $X_t^v$ until the rates of the two processes differ) do not coincide on $t \in [0,1]~$.}
{By taking $\tau$ as the identity map in \eref{e:dj},} we further have for all $v$ that
$$
d_{J_1}(\widetilde X^v,X^v) \le
\sup_{t \in [0,1]} \|{\widetilde X_t^v - X_t^v }\|_1 \le v^{-1} {\bar c_\bullet}
$$
\noindent
and consequently the required $J_1$-weak \abbr{ldp} for {$\{\widetilde{X}^v_t\}$}
follows from the \emph{local} \abbr{ldp}
for {$\{ X^{v}_t \}$} with respect to the $d_{J_1}$-metric balls
(see \cite[Thm.~4.1.11]{LDT1}). In view of
\eref{e:exp-good-apx}, the latter local \abbr{ldp} follows from having
for any $z \in D_{0,1}(\RR^d)$ and
{all $\rho$ large enough (which may depend on $z(\cdot)$),}
\begin{align}\label{t:loc-ubd}
\inf_{\delta>0}
\limsup_{v \rightarrow \infty} \frac{1}{v} \log \px{x_0^v}{d_{J_1}(
X^{v,\rho},z) < \delta}
& \le - I_{x_0,1} (z) ~\,,
\\
\inf_{\delta>0}
\liminf_{v \rightarrow \infty} \frac{1}{v} \log \px{x_0^v}{d_{J_1}(X^{v,\rho},z) < \delta}
& \ge  - I_{x_0,1}(z) \,.
\label{t:loc-lbd}
\end{align}
In establishing these bounds we tackle the diminishing rates
$\lambda_r(\cdot)$ at $\partial \RR^d$ by employing
a \abbr{ldp} from \cite{dupuis16} for the empirical measure sample-path
$t \mapsto \mu^n_t$ of $n$ mean-field interacting particles. Specifically,
fixing
$z \in D_{0,1}(\RR^d)$ let $\gamma := 1 + \sup_{t \in [0,1]} \|z(t)\|_1$.
Since for any $v$ and $\rho$
\begin{equs}\label{e:gamma-enough}
 d_{J_1}(X^{v,\rho},z)<1  \qquad \Longrightarrow \qquad
\{ X_t^{v,\rho} : t \in [0,1] \} \subseteq \widetilde K_\gamma \,,
\end{equs}
the choice of jump rates $\Lambda^{v,\rho}(\cdot)$
outside
$\widetilde K_\gamma$ \emph{is irrelevant} for the bounds
\eref{t:loc-ubd} and \eref{t:loc-lbd}. Choosing an
integer $\rho$ with $\rho \ge 2 \rho_0(\gamma) \ge 2 \gamma$ which is further
large enough for
\eref{dfn:rho-val} to hold, the process $X_\cdot^{v,\rho}$
is confined to $\widetilde{K}_\rho$ of \eref{dfn:K-tilde}
so has at most $n=v \rho$ molecules
(to simplify notations, take \abbr{wlog} $v \in \mathbb N$).
We thus consider the evolution of $n$ indistinguishable particles,
each labeled by a type from $\S_\star$, where $n \mu^n_t(\star)$
counts the $\star$-particles that compensate the
$c^r_\star$ molecules gained/lost at each reaction.
{Starting at $v (x^{v}_0)_i$ particles of type {$s_i \in \S$} and
$n - v \|x^v_0\|_1$ of $\star$-type,}
our goal is to have for $\Psi(\cdot)$ of \eref{e:bijec-ode}
the continuous bijection
\begin{equs}\label{e:bijec}
X_{t}^{v,\rho} = \Psi(\mu^n_t) \,, \qquad
\mu^n_t (\star) = 1 - \rho^{-1} \| X_{t}^{v,\rho} \|_1 \,.
\end{equs}
To this end, a chemical reaction $r \in \R$ is mapped to the simultaneous change
of $\ell_r := \|c^r_\inn\|_1 \vee \|c^r_{\rm out}\|_1 \le { {\bar c_\bullet}}$ particle types, where, given $\mu^n$, any ordered $\ell_r$-tuple
${\bf i} \in \S_\star^{\ell_r}$ that has type-count
configuration $((c^r_\star)_+,c^r_\inn)$ independently changes
into an ordered $\ell_r$-tuple ${\bf j} \in \S_\star^{\ell_r}$
that has type-count configuration $((c^r_\star)_-,c^r_{\rm out})$, at the rate
\begin{equs}\label{e:gamma-rate}
\Gamma^{(r),n}_{{\bf i}{\bf j}}(\mu^{n})
= \frac{k_r \, \ell_r! \, v^{1-\|c^r_\inn\|_1}}
{M_r
\binom {n \mu^n(\star)}{ (c_\star^r)_+}
 (c_\star^r)_+ !}
\,,
\end{equs}
where $M_r = \ell_r!^2 / [(|c^r_\star|)! \prod_{i=1}^d (c^r_\inn)_i! (c^r_{\rm out})_i!]$ is the number of pairs $({\bf i},{\bf j})$ matching
the specified type-count configurations (and to accommodate all possible
\abbr{crn}s we permit ${\bf i}_l={\bf j}_l$ for some $l$, unlike
\cite[Eqn. (2.1)]{dupuis16}).
Indeed, for
$\Lambda_r^{v,\rho}$ of \eref{e:dfn-Lambda-rho} and
$\{\Gamma^{(r),n}_{{\bf i}{\bf j}}\}$ of \eref{e:gamma-rate},
the generator of $\mu^n$ in \cite[Eqn. (2.7)]{dupuis16} has total
jump rate $v \Lambda^{v,\rho}_r ( \Psi(\cdot) )$
in each direction {$(-c^r_\star,c^r)$}, $r \in \R$, thereby
yielding the bijection property \eref{e:bijec}.
From \eref{e:gamma-rate} it is also
easy to check that for any $\mu^n \to \mu$,
\begin{equ}
n^{\ell_r-1} \Gamma^{(r),n}_{{\bf i}{\bf j}}(\mu^n) \to
\tilde k_r \mu_\star^{{-(c^r_\star)_+}}
=: \Gamma^{(r)}_{{\bf i}{\bf j}}(\mu)
\,,
\end{equ}
where $\tilde k_r > 0$ is independent of $\mu$.
Such $\{\Gamma^{(r)}_{{\bf i}{\bf j}} (\mu) \}$ satisfy
the uniformity condition of \cite[Assmp.~3.1]{dupuis16}.
{On
$\mathcal M_+ (\S_\star) := \{ \mu \in \mathcal M_1(\S_\star): \mu_\star \ge 1/2 \}$
they also have the Lipschitz continuity of \cite[Assmp.~2.2]{dupuis16},
and taking into account the factor $v/n$ between volume normalizations,
we have on $\mathcal M_+(\S_\star)$ the
Lipschitz continuous asymptotic normalized reaction rates
$\rho^{-1} \lambda_r(\Psi(\mu))$ for $\mu^n$ that
satisfy \cite[Property 2.3]{dupuis16}.} As shown in
\cite[Section 6]{dupuis16}, having \cite[Property 2.3]{dupuis16}
throughout $\mathcal M_1(\S_\star)$ yields
the \abbr{ldp} upper bound for $\{\mu^n_t\}$ in the $J_1$-topology of
$D_{0,1}(\mathcal M_1(\S_\star))$, at rate $n$.
{Here $\mu^n_0 \to \Psi^{-1}(x_0)$ and} the asymptotic
reaction rates for $\mu^n$ depend only on $\Psi(\mu)$.
Consequently, the rate function controlling the \abbr{ldp} upper bound
for $\{\mu^n(t)\}$ is
\[
J(\mu) = \rho^{-1} I_{x_0,1}(\Psi(\mu)) \,,
\]
and upon compensating for the factor $v/n$ between the two rates,
such an \abbr{ldp} upper bound for $\{ \mu^n(t) \}$ readily yields
\eqref{t:loc-ubd}.
Our problem fails to satisfy the Lipschitz continuity of
\cite[Property 2.3]{dupuis16} when $\mu_\star=0$. However,
$\rho \ge 2 \gamma$ guarantees that $\mu_\star \ge 1/2$
on $\Psi^{-1}(\widetilde K_\gamma)$, which in view of
\eqref{e:gamma-enough} is all that matters for
\eqref{t:loc-ubd}.
{As explained at the start of the proof, upon combining \eqref{t:loc-ubd}
with the exponential tightness of \lref{l:new},} {we get the stated \abbr{ldp} upper bound of \eref{dfn:ubd}} {(for $T=1$). In particular,
due to exponential tightness the \abbr{lhs} of \eqref{dfn:ubd} is
zero for some compact $\Gamma$. The same applies then for the infimum of
$I_{x_0,1}(\cdot)$ over this compact set}
{and hence $I_{x_0,1}(\zeta) =0$ for some
$\zeta \in AC_{0,1}(\RR^d)$ with $\zeta(0) = x_0$. Recall
\rref{r:zero-rate} that such $\zeta(\cdot)$ must satisfy
the \abbr{ode} \eref{e:ma} for $t \in [0,1]$. We note in passing that
the same argument applies for any finite $T$ (as explained just
prior to \lref{l:new}), yielding the existence of
global solutions for this \abbr{ode}
and further, the \abbr{flln}
of \rref{r:boundednessoftrajectories} then holds when
$v_k/(\log k) \uparrow \infty$, by combining the
stated \abbr{ldp} upper bound and the Borel-Cantelli lemma.}

Next, note that
we have \eqref{t:loc-lbd} as a consequence of the \abbr{ldp} lower
bound of \cite[Eqn. (8.1)]{dupuis16} holding for {$\{\mu^n_t\}$}.
As mentioned in \cite[Remark 8.6]{dupuis16}, such \abbr{ldp} applies
when having in addition to \cite[Property 2.3 \& Assmp.~3.1]{dupuis16}
also the $\eta$-ergodicity of \cite[Assmp.~3.3]{dupuis16}
and that the solution of the \abbr{ODE} \eqref{e:lln-dupuis}
satisfies \cite[Property 4.13]{dupuis16}. The latter amounts
to having the lower bound of \eqref{e:bd-413} also
for $\mu_\star(t)$. Having $\rho \ge 2 \rho_0(\gamma)$, from
\lref{l:p314dupuis} this holds whenever starting at
$\mu(0)\vert_{\S} \le \gamma/\rho$ which is precisely
$\Psi^{-1}(\widetilde K_\gamma)$ (and thus all that
is relevant for
\eqref{t:loc-lbd}).
The $\eta$-ergodicity of \cite[Assmp.~3.3]{dupuis16} amounts here
to being able to reach a particle population that exhibits all $d+1$
elements of $\S_\star$ upon starting at $n \gg 1$ particles
from a fixed, single type from $\S_\star$ and
\aref{a:1}(b) guarantees this when starting at only $\star$-particles.
We thus have also \cite[Assmp.~3.3]{dupuis16} except at the face
$\mu_\star=0$ on the boundary of $\mathcal M_1(\S_\star)$. While the
behavior at that face plays a role for some events, thanks to
\eqref{e:gamma-enough}
it is irrelevant here.
\end{proof}

\section[The stability of {\lowercase{\sc{ASE}}} networks]{The stability of {\sc{\bf ASE}} networks}
\label{s:prop}

{The proof of \pref{t:thm}
is long and technically challenging, so we first
sketch in Section \ref{s:prop2} the proof
of \eqref{e:expdrift} for $x$ away from $\partial \RR^d$
to familiarize the reader with the techniques used in the subsequent sections, where we carry out this proof in full detail.
}

\subsection{Toric rays and outline of the stability proof}\label{s:prop2}

Following the geometrical analysis of \cite{gopal13}, we
first define \emph{toric rays}, using throughout
for $w \in \mathbb R^n$, $z \in (\RR^n)^o$ and $\theta\in\RR$ the operators
\begin{align*}
\log(z) &:= (\log z_1, \dots , \log z_n) \in \mathbb R^n~,\\
z^w &:= (z_1^{w_1},\dots ,z_n^{w_n}) \in (\mathbb R_{+}^n)^o~,\\
\theta^w &:= (\theta^{w_1},\dots, \theta^{w_n}) \in (\mathbb R_{+}^n)^o~.
\end{align*}
\dfn{To each $w$ in the unit sphere $S^{n-1}$ we associate the $w$-\emx{toric ray}
$$
T^w = \bigcup_{{\theta \geq 1}}\,\theta^w \subset \RR^n \,.
$$
We also introduce the toric-ray parameters
\begin{align}
 \theta(z) := \exp(\|\log(z)\|_2)\,, \qquad
w(z) :&=\frac{1}{\log \theta(z)} \log(z)\,, \nonumber \\
\pc{\theta,w} : \pc{\RR^{n}}^o \rightarrow  \Rr_{>1} \times S^{n-1}~,
\quad z &= \theta(z)^{w(z)}~\,.
\label{dfn:topic-p}
\end{align}
\Label{d:toric-ray}
}

\rmk{To see why $U(\cdot)$ of \eqref{e:lyapunov1} is most suitable
for mass action systems, note that along a $w$-toric ray
\begin{equ}
\nabla U(\theta^w) = \log(\theta^w) = (\log \theta) w \,,
\label{e:grad-U}
\end{equ}
while the derivative of the \abbr{ode} \eqref{e:ma} at a point on such
a ray is
\begin{equ}
\frac{\d \zeta}{\dt} \Big|_{\zeta = \theta^w} = \sum_{r \in \R}
\lambda_r(x) c^r \Big|_{\zeta = \theta^w} = \sum_{r \in \R} k_r
\left(\theta^w\right)^{c^r_{\rm in}} c^r = \sum_{r \in \R} k_r \theta^{\langle w , c^r_{\rm in} \rangle} c^r~\,.
\label{e:asymptotic}
\end{equ}
Thus,
at $x = \theta^w$
the time derivative of $U(\xt )$ for the solution $\xt $ of \eqref{e:ma}
is
\begin{equ}
\frac{\d}{\d t}U(\xt )\Big|_{\zeta = \theta^w} = \langle \nabla U(x) , \frac{\d \zeta}{\d t} \rangle \Big|_{\zeta = \theta^w}= (\log \theta) \sum_{r \in \R} k_r  \langle w , c^r \rangle \theta^{\langle w , c^r_\inn \rangle}~\,.
\label{e:lyapunovdominance}
\end{equ}
For fixed $w$ and $\theta \gg 1$ the sum on the \abbr{rhs} of
\eref{e:lyapunovdominance} is dominated by reactions $r \in \R_{w}$
{(maximizing $\langle w , c^r_\inn \rangle$).} Thus, in
strongly endotactic \abbr{crn}s, where at least one such reaction
contributes negatively to this sum
by having $\langle w , c^r\rangle < 0$, and no other reaction $r$ in
$\R_w$ contributes positively to it, the \abbr{lhs} of
\eref{e:lyapunovdominance} will also be negative for all large enough $\theta$.
As shown in \cite{gopal13}, if this applies uniformly over
$w \in S^{d-1}$ then for some compact $K$ we have
$\frac{\d}{\d t} U(\xt ) < 0$ whenever $\xt  \notin K$, so
\eqref{e:ma} has {an absorbing compact set}.
Indeed, suppose to the contrary, that for some diverging sequence $x(j) \in \RR^d$
\begin{equ}
\frac{\d}{\dt} U(\xt )\Big|_{\zeta = x(j)} \geq 0 \qquad \forall\, j \in \mathbb N~\,. \label{e:contra}\end{equ}
By compactness of $S^{d-1}$, upon passing to a suitable sub-sequence, the
corresponding toric-ray parameters $x(j)=\theta(j)^{w(j)}$
{form a \emph{toric-jet} of frame $\bar w
 = \{ \bar w^{(k)} : k \le \ell\}$
(see \dref{dfn:toric-jet} and \cite[Lemma 6.7]{gopal13}),
where $w(j) \to {\bar w}^{(1)}$ and $\theta(j) \to \infty$. By compactness of
$[1,\infty]$, there exists a further sub-sequence along which
$x(j)^{\bar w^{(k)}}$ converge for each $k$
(possibly to $\infty$), implying the convergence of the functions $\widehat \varphi_{r}(x):=
k_{r} \langle w,c^r \rangle \theta^{\langle w , c^{r}_\inn \rangle}$.}
%
For strongly endotactic \abbr{crn}s one
can show \cite{gopal13} that along
such a
toric-jet, for any $r \in \R$ there exists $r' \in \R$ whose contribution $\widehat \varphi_{r'}(x(j))$
to the \abbr{rhs} of \eqref{e:lyapunovdominance} is such that
$\lim_j \widehat \varphi_{r'}(x(j)) < 0$ and
{$-\widehat \varphi_{r'}(x(j))/(\widehat \varphi_{r}(x(j)))_+ \to \infty$ (where $0^{-1}:= \infty$),}
contradicting \eref{e:contra}.
\label{r:lyapunovdominance}}

{\rmk{Note that in the components $s_i \in \S$ where $w_i$ is negative, the value of $\theta^{w_i}$ decreases as $\theta \to \infty$. Such $w$ are therefore used to parametrize through \eref{dfn:topic-p} points in $(v^{-1} \mathbb N_0)^d$ that are at a distance $< 1$ from the boundary $\{x_i = 0\}$.}}

\noindent
\pref{t:thm} amounts to having for some finite $b$, for any
$x \in \Ov$ and for $v>\vp (\|x\|_1)$,
\begin{equ}
\sum_{r \in \mathcal R} \Lambda_r^{(v)}(x)
\Big[ U^v(x+v^{-1} c^r) - U^v(x) \Big] \leq  e^{bv}\,.
\label{e:contra2}
\end{equ}
Recall that $\Lambda_r^{(v)}(x) \le \lambda_r(x)$ which is uniformly bounded
on compacts, as is $U(x)$. Hence, there exists a finite $b=b(\rho)$  such
that \eqref{e:contra2} holds for any $v \ge 1$ whenever $\|x\|_1 \le \rho$.
Letting
\begin{equa}
{\AN^v} &:= {\{ x \in \Ov : \rho < \|x\|_1 \le \rho' \} \,,}\Label{dfn:anv}\\
L_{r}^{(v)}(x) &:= U(x) (Q_{r}^{(v)}(x)-1)\,, \qquad
Q_r^{(v)} := U^v(x+v^{-1}c^r)/U^v(x) \,,\Label{dfn:L-and-Q}
\end{equa}
we thus establish \pref{t:thm} upon showing that for some $\rho<\infty$
any $\rho' \ge \rho$, $x \in \AN^v$ and $v > \vp (\rho')$, we have
\begin{equs}\label{dfn:avx}
a^{(v)} (x) : = \sum_{r \in \R} k_r \varphi^{(v)}_r(x) \le 0 \,, {\quad
\varphi^{(v)}_r(x):= k_r^{-1}  \Lambda_r^{(v)}(x) L_r^{(v)}(x) }\,,
\end{equs}
where, by \eqref{e:jrates} one considers in $a^{(v)}(x)$ only $r$
such that $v x \ge c^r_\inn$ (thus $x + v^{-1} c^r \in \RR^d$).
Subject to having the $v$-independent approximation for all $x \in (v^{-1} \mathbb N)^d$,
\begin{equ} Q_r^{(v)}(x)
= \exp\Big[ \frac{h_r(x) + \mathcal O_{\|x\|}(1)}{U(x)} \Big] \qquad \text{with} \qquad h_r(x):=\langle \nabla U(x), c^r \rangle\,,
\label{dfn:Qrv}
\end{equ}
we can prove \eqref{dfn:avx},
{at least for a strictly positive $x$},
by contradiction. Specifically, {one can show that
it suffices to rule out having
$a^{(v(j))}(x(j)) > 0$ along a \emph{rapidly diverging} volume-jet
$(v(j),x(j))$. That is, along some diverging toric-jet
$x(j)=\theta(j)^{w(j)} \in (v(j)^{-1} \mathbb N)^d$,
with $\theta(j) \to \infty$ and frame $\bar w$,
such that $v(j) \to \infty$ arbitrarily fast (\ie
allowing for an arbitrary $\vp (\rho)$ in \dref{d:vjss}).}
Similarly to
\rref{r:lyapunovdominance},
we arrive at a contradiction by showing
that for {some $\vp $ any such $\vp $-divergent volume-jet $(v,x)$
and} $r \in \R$, there must exist
$r' \in \R_{ {\bar w^{(1)}}}$ such that
eventually $\varphi^{(v)}_{r'}(x) < 0$ and
$-\varphi^{(v)}_{r'} (x)/(\varphi^{(v)}_r(x))_+ \to \infty$.
To this end, we first show in \lref{r:noindicatoratboundary} that
{for $\vp (\rho)=e^{\rho}$ and some constants $\delta_{r'}>0$,
along any $\vp $-divergent volume-jet $(v,x)$ framed by $\bar w$, eventually
\begin{equ}
\Lambda_{r'}^{(v)}(x) \ge \delta_{r'} \lambda_{r'}(x)
 \qquad \forall\,r' \in \R_{\bar w^{(1)}} ~,
 \label{e:lambdaw}
\end{equ}
}
which as $\Lambda_r^{(v)}(x) \leq \lambda_r(x)$, implies that for
$C<\infty$, any $r \in \R$ and $r' \in \R_{\bar w^{(1)}}$,
{eventually}
\begin{equs}\label{dfn:Prr}
C \pd{\frac{ \varphi^{(v)}_{r'}(x)}{{\varphi^{(v)}_r(x)}}}
\geq
\frac{k_{r'}^{-1}\lambda_{r'}(x)}{k_r^{-1}\lambda_r(x)}
\pd{\frac{L_{r'}^{(v)}(x)}{L_{r}^{(v)}(x)}}
=: P^{(v)}_{r,r'}(x) \,.
\end{equs}
Referring to the first part in the \abbr{rhs} of \eref{dfn:Prr} as a
\emx{monomial term} (since $k_{r'}^{-1}\lambda_{r'}(x)/k_r^{-1}\lambda_{r}(x)
\\ = \theta^{\langle w, c_\inn^{r'}-c_\inn^r\rangle}$),
and to the second part (in the absolute value sign) as \emx{Lyapunov term},
we then show that for any
$r \in \R$, if eventually $L_r^{(v)}(x) > 0$ then
by \cite[Prop.~6.24]{gopal13} there exists {$r' \in \R_{\bar w^{(1)}}$}
with $h_{r'}(x) \rightarrow -\infty$ such that along the
{divergent volume-jet,}
\begin{equs}\label{dfn:Prrp}
\lim_{j \to \infty} P^{(v(j))}_{r,r'}(x(j)) = \infty \,.
\end{equs}
Indeed, relying on \eref{dfn:Qrv} we establish \eqref{dfn:Prrp}
by proceeding according to whether
$\kappa_r := \lim_{j} {| h_r(x{{{}}}) |}$ is finite
or infinite.
Specifically, we have the following cases:
\begin{myenum}
\item Lyapunov domination, where $\kappa_r$ is finite and with
$U(x) \to \infty$ we have that
$L_r^{(v)}(x) $ remains uniformly bounded.
The existence of $r'(r) \in \R_{\bar w^{(1)}}$ with $L_{{r'}}^{(v)}(x) \rightarrow -\infty$ resulting from \cite[Prop.~6.24]{gopal13} (see \lref{l:innerproductdomination}), combined with \cite[Lemma 6.22]{gopal13} to bound
the monomial term away from zero, concludes the proof.
\item Monomial domination, where $\kappa_r = \infty$ so
$Q_r^{(v)}(x{{{}}}) = e^{h_r(x{{{}}}) /U(x{{{}}})}(1+o(1))$ by \eref{dfn:Qrv}. This implies, by \cite[Prop.~6.20 \& 6.24]{gopal13}, the existence of $r' \in \R_{\bar w^{(1)}}$ such that $|L_{r'}^{(v)}(x)/L_{{r}}^{(v)}(x)| = {\mathcal O} (\theta^{-\langle w, c^{r} \rangle/U(x)})$, whose exponent goes to zero as $j\rightarrow \infty$. On the other hand, for such $r'$ by \cite[Lemma 6.22]{gopal13} the exponent of $\theta$ in the monomial term of \eref{dfn:Prr} is (eventually) strictly positive and bounded away from zero along the
{toric jet,} thereby establishing \eref{dfn:Prrp}.
\end{myenum}

In order to establish \eqref{dfn:avx} also on $\partial \RR^d$, we
need to extend the preceding program to deal with boundary effects such as the divergence of
$\nabla U(x)$. This is done by separately considering each
face of $\RR^d$.
In particular, Section~\ref{s:approximations}
establishes
\eref{dfn:Qrv}
in a more general form,
substituting $\nabla U(\cdot)$
{with the \emph{$v$-dependent}}
$\nabla_r^{(v)}U(\cdot)$
of \eqref{dfn:gradU-P}.
Section~\ref{s:gopal} adapts the definitions of toric jet and strongly endotactic \abbr{crn}s from \cite{gopal13}
as needed for $\partial \RR^d$. This and the corresponding
results from \cite[Sec. 6]{gopal13} are then used in
Section~\ref{s:tj} to show the divergence of {$\prrr$, first for
Lyapunov domination (in \lref{l:innerproductdomination3}),
and then for monomial domination (in \lref{l:monomialdomination}). Finally,
Section~\ref{s:proofofprop} follows the preceding outline in combining
everything to a proof of \pref{t:thm}.}

\subsection{Approximation lemmas}
\label{s:approximations}

The image of $\Rr^d$ under the exponential map is $(\RR^d)^o$, so we
will establish \eqref{dfn:avx} separately on the various faces of $\partial \RR^d$
by considering the relevant
\abbr{crn}s {$(\S,\C,\R(\PP))$ where, for any non-empty $\PP \subseteq \S$,
\[
\R(\PP) := \{ r \in \R ~:~ \supp\{ c^r_{\inn}\} \subseteq \PP\}
\]
denotes the reactions with substrates only from
 $\PP$.}
To this end, identify such $\PP = (s_{i_1},\dots,s_{i_{d_\PP}})$
of cardinality $|\PP| = d_\PP {\ge 1}$
with the corresponding indices $(i_1,\dots,i_{{d_\PP}})$, denoting by
$\mathbb S^d(\PP)$ the restriction of a space $\mathbb S^d$
(be it $\Rr^d$, $\RR^d$, $\mathbb N_0^d$ or $\mathbb N^d$),
to these coordinates (\ie having zero values outside $\PP$).
Aiming at the approximation \eqref{dfn:Qrv}
for {$x \in (v^{-1} \mathbb N)^d(\PP)$} and $r \in \R(\PP)$, we modify $\nabla U(x)$ to
\begin{equs}\label{dfn:gradU-P}
\pc{\nablavU}_i :=
\left\{\begin{array}{c@{,}l}
		\log x_i	         & 	\qquad i \in \PP  \,	\\
		\log (v^{-1} c^r_i)  & \qquad i \in \supp \{c^r_\outt \} \cap \PP^c \,\\
		0                    &  \qquad \text{ otherwise }
	\end{array}\right..
\end{equs}
We write $\epsilon_v(x)$ for functions that are uniformly bounded
in $x$ by some $\bar \epsilon_v \to 0$ as $v \rightarrow \infty$
and $\epsilon(x)$ for any globally bounded function of $x$.
\lmm{Setting $g_p(x):=\|\log(x)\|_p$\, for $p=1,2$, we have that
\begin{equ}
\frac{2 g_1(x)}{v U(x)} \le
\frac{d + g_2(x)^2}{v U(x)} \le \epsilon_v(x) \qquad \forall x \in \Ovp ~.
\label{e:lexpisok}
\end{equ}
\label{l:expisok}}
\begin{proof} {The first inequality in \eref{e:lexpisok} directly results from the fact that $x^2+1\geq2|x|$ for all $x \in \Rr$.} {Next},
since $g_2(x) \le \sqrt{d} \sup_i \{|\log x_i|\}$
and $U(x) \ge 1$, by \eqref{e:lyapunov1} it suffices
to show that
\[
\frac{|\log y|^2}{v [y(\log y -1) + 2]} \le \epsilon_v(y) \qquad
\qquad \forall y \ge v^{-1} \,.
\]
For $y \in [v^{-1},v]$ the \abbr{lhs} is at most
$(\log v)^2/v \to 0$
as $v \to \infty$, whereas for $y \ge v \ge e^2$ the \abbr{lhs} is bounded
above by $2\log y/(vy) \le 2 \log v/v^2 \to 0$ as $v \to \infty$.
\end{proof}

\rmk{Hereafter, {for any $r \in \R(\PP)$,} we consider \abbr{wlog} only {$r$-\emph{relevant}} $x$,
namely those for which
$vx + c^r \in \mathbb N_0^d$, for otherwise the
corresponding term disappears in \eqref{dfn:avx} (see
\rref{r:zerosummands}).}

\lmm{There is a finite $v_*$ such that for any $\PP\subseteq \S$, all
  $r \in \R(\PP)$ and all {$r$-relevant} $x \in \Ovp(\PP)$, one has
  for $v\ge v_*$:
\begin{align*}
{Q_r^{(v)}(x)}
& = \exp \pq{\frac{ \h{v}{x}
+ \epsilon(x)}{U(x)}}\,,
\end{align*}
where $\h{v}{x} := \langle \nablavU, c^r \rangle$\,.
\label{l:taylor}
}

%

\begin{proof} Since the number of possible $\PP$ and $r$ is finite, it
  suffices to prove the claim for fixed $\PP$ and $r$. We have in terms of $f:= v[ U(x+ v^{-1} c^r) - U(x)]/U(x)$ that
\begin{equ}
{Q_r^{(v)}(x)} = \Big(1+ \frac{f}{v}\Big)^v = \exp\Big[ f - v  R(f/v) \Big]
\end{equ}
where the non-negative $R (y) := y - \log (1+y)$ satisfies
\begin{equs}\label{e:R-bd}
 R(y) \le 2 y^2 \,,
\qquad \forall y \ge -1/2 \,.
\end{equs}
Now, for any $r \in \R(\PP)$ and $x \in \Ovp(\PP)$ with
$vx + c^r \in \mathbb N_0^d$ we have that
\begin{equ}
f \, U(x) - \h{v}{x}
= \sum_{i \in \PP} \psi(v x_i;c^r_i)
- \langle c^r, \onevec \rangle = \epsilon(x)
\label{e:rU/U=log}
\end{equ}
is uniformly bounded
since $\psi(b;c):= (b+c)\log(1+c/b)$ decreases in $b \ge \max(1,-c)$.
Hence, from \eqref{e:rU/U=log}, \eqref{dfn:gradU-P} and \lref{l:expisok},
\[
\frac{1}{2} f^2 \le \Big(\frac{\h{v}{x}}{U(x)}\Big)^2
 + \Big(\frac{\epsilon(x)}{U(x)}\Big)^2
= \frac{v \epsilon_v(x)}{U(x)} \,.
\]
Since $U(x) \ge 1$ we see that
$(f/v)^2 \le 2\epsilon_v(x)/v \le 1/4$ for
some $v_*$ finite, any $v \ge v_*$ and all $x$,
in which case by \eref{e:R-bd} we have that
$v R(f/v) \le 2 f^2/v
\le 4\frac{\epsilon_v(x)}{U(x)}$, as claimed.
\end{proof}

\subsection{Strongly endotactic sub-networks and divergent volume-jets}
\label{s:gopal}
Throughout, for {non-empty} $\PP \subseteq \S$ and $w \in\Rr^d$ we denote by $\pi_\PP: \Rr^d \to \Rr^{d_\PP}$ the projection onto the coordinates with indices in $\PP$.
Proceeding to adapt for $(\S,\C,\R(\PP))$ key definitions
from \abbr{crn} theory, such as \emph{strongly endotactic} (see \cite{gopal13}),
for all $w\in\Rr^d$ with non-zero projection $w_\PP := {\pi_\PP} w$,
let $\R(\PP)_w$ denote the reactions having $c^r_\inn$
in the $w$-maximal subset of
$\C_\inn(\PP) = \{c^r_\inn : r \in \R(\PP)\}$.
Clearly, $\R(\PP)_w$ depends only on $w_\PP$
which \abbr{wlog} is in the
$(d_\PP-1)$-dimensional unit sphere $S^{\PP}$
and \abbr{wlog} we further identify $\C_\inn(\PP)$
with ${\pi_\PP} \C_\inn(\PP)$.
\dfn{Fixing $w_\PP \in S^{\PP}$, a reaction $r \in \R(\PP)$ with
$\supp\{ c^r_{\rm out} \} \subseteq \PP$ is called
\emph{$w$-dissipative}, \emph{$w$-null} or \emph{$w$-explosive}
according to $\langle w, c^r \rangle = \langle w_\PP, {\pi_\PP} c^r \rangle$ being
negative, zero or positive, respectively. Any
$r \in \R(\PP)$ having some product species not within $\PP$
is considered \emph{$w$-dissipative} (regardless of $w$).
Similarly, $r \in \R(\PP)$ is $\{w\}$-explosive,
$\{w\}$-null or $\{w\}$-dissipative, if the relevant
property holds for all but finitely many elements
of $\{w\} \subset S^{\PP}$.
\label{d:dissipative}}
\begin{remark}
For $\PP=\S$ our \dref{d:dissipative} of $w$-dissipative and $w$-explosive
reactions, coincides with \cite[Def. 6.15]{gopal13} of
$w$-sustaining and $w$-draining reactions, respectively.
The nomenclature was changed to stress the behavior of reactions
for $\|x\|_1 \gg 1$ which is of interest here: dissipative [explosive] reactions contribute to the decrease [increase] of the Lyapunov function along trajectories far away from the origin.
\end{remark}
We next extend  \dref{d:Pendo1} of strongly $\S$-endotactic \abbr{crn}
to $\PP \subset \S$. Such an extension is needed in light of
\rref{r:lyapunovdominance}, and made relevant by \lref{l:endotacticorpositive}.
\dfn{For any $w\in \Rr^d$ with non-zero projection onto $\PP$ (or $w_\PP \in S^\PP$), the \abbr{crn} $(\S,\C,\R(\PP))$ is called
\emph{$w$-strongly $\PP$-endotactic} if
the set $\R(\PP)_{w}$ contains at least one $w$-dissipative reaction,
and no $w$-explosive reactions. Such a \abbr{crn}
is called \emph{strongly $\PP$-endotactic} if it is $w$-strongly
$\PP$-endotactic for all  $w$ as above.
\label{d:Pendo}}
\lmm{\label{l:endotacticorpositive}
Any strongly endotactic \net is strongly $\PP$-endotactic for
all $\PP \subset \S$ {if $\R(\PP) \ne \emptyset$}.}
\begin{proof}
{Fixing $\PP \subset \S$ with $\R(\PP) \ne \emptyset$},
suppose that for a non-zero $w$ there is a $w$-explosive 
$r \in \R(\PP)_w$. Since the non-negative
$\sum_{i \notin \PP} (c^{r'}_\inn)_i$ is zero iff $r' \in \R(\PP)$,
setting $w'_i = w_i$ for $i \in \PP$ and
$w'_i = - \gamma$ for $i \notin \PP$ we have that
$\R_{w'} = \R(\PP)_w$ for $\gamma$ large enough.
Further, $\supp\{ c^{r} \} \subseteq \PP$ hence
$\langle c^{r}, w' \rangle = \langle c^{r}, w \rangle > 0$,
so having $r \in \R_{w'}$ contradicts \dref{d:Pendo1}.
For the same reason, if every reaction in $\R(\PP)_w$ is $w$-null,
then for large $\gamma$ the same applies for every reaction in $\R_{w'}$,
in contradiction with \dref{d:Pendo1}.
\end{proof}
To show that \eref{dfn:avx} holds {\emph{whenever}
$v > \vp (\|x\|_1)$ and $v x \in {\mathbb N}^d(\PP)$ with
$\|x\|_1 \ge \rho$,}
requires an approximation framework
for sequences $x(j) = \theta(j)^{w(j)}$ satisfying
$\theta(j) \to \infty$ and $w(j) \to \bar w^{(1)}$ in $S^\PP$.
To this end, we follow \cite[Sec. 6]{gopal13} in coding
the latter convergence by a suitable $\dsr$-dimensional \emph{frame} \cite{miller09}: an
orthonormal system (\abbr{ons})
$\bar w := \{\bar w^{(1)},\dots ,\bar w^{(\dsr)} \} \subset S^{\PP}$
such that
\begin{equs}
\lim_{j \rightarrow \infty} \frac{\beta^{(k+1)}{{{}}}}{\beta^{(k)}{{{}}}}
= 0 \,, \qquad \forall k < \dsr \,, \quad
\beta^{(k)}{{{}}} = \beta^{(k)}(j) := \langle w(j), \bar w^{(k)}\rangle \,.
\label{dfn:u-jet}
\end{equs}
For generic $\{w(j)\}$ one needs a full $d_\PP$-dimensional basis
of $S^\PP$, but degeneracy allows for $\dsr < d_\PP$ (\eg,
$\bar w^{(1)}$ alone suffices when all $w(j)$ lie on a single
toric-ray). Further, the \emph{order} within $\bar w$
is adapted to the sequence, so that the angle between $w(j)$
and $\bar w^{(k)}$ decays faster with each increase of $k$. {Through the following definition, by slight abuse of notation we suppress the index $j$ for elements of the sequence $\{x(j)\}$ and other related quantities to increase the readability of forthcoming formulas.}
\dfn{{\rm \cite[Defn. 6.2, 6.18]{gopal13}}~
\newline
(a) A \emx{unit jet} on {a frame} $\bar w$ is a sequence $\{w\} = \{w(j)\}$ of
unit vectors in the {conic hull} $\mathcal {\rm Co}(\bar w)$ satisfying
\eqref{dfn:u-jet}.
\newline
(b) A \emx{toric jet} $\{x\}$ is a sequence
$\theta{(j)}^{w(j)} \in \Rr_{>0}^d(\PP)$ for a unit jet
$\{w{{{}}}\}$ and $\theta(j) \to \infty$\,.
\newline
(c) A unit jet $\{w{{{}}}\}$ and the corresponding toric jets
are \emph{adapted} to $(\S,\C,\R(\PP))$ if the classification
of each $r \in \R(\PP)$ according to \dref{d:dissipative} is
conserved by all {$w(j)$ with} $j \in \mathbb N$ and for all $k = 1, \dots, \dsr$, $\lim_j \theta{{{}}}^{\beta^{(k)}{{{}}}}$ exists and takes values in $[1,\infty]$~.
\label{dfn:toric-jet}
}

\begin{remark}
When the unit jet $\{w\}$ is adapted to $(\S,\C,\R(\PP))$ and
clear from the context, in view of point (c) we call a reaction
$r \in \R(\PP)$ dissipative or explosive, per \dref{d:dissipative}, without explicitly indicating the choice of $w(j)$.
\end{remark}

Having assigned any $r \in \R(\PP)$ with $\supp \{c_\outt^r\} \not\subseteq \PP$
as dissipative reactions, it is necessary for the
strategy presented in Sec.~\ref{s:prop2} to ensure that their contribution to $a^{(v)}(x)$ is negative along $\{(v,x)\}$ in case
$r \in \R(\PP)_{\bar w^{(1)}}$.
Since {for such a reaction}
$\lim_v \langle \nabla_r^{(v)} U(x), c^r \rangle
= - \infty$, we obtain in \lref{l:e-1/22}
the desired behavior of $L_r^{(v)}(x)$
by choosing, for every $x$ to have $v > \vp (\|x\|_1)$ for a
function $\vp (\cdot)$ increasing fast enough. Our next definition
guarantees that this condition on $v$ is met along $\{x\}$.
\dfn{Fixing $\PP \subseteq \S$ and an increasing function $\vp (\rho)$, we
call a sequence of $\{(v,x) : v x \in {\mathbb N}^d(\PP),
{v>\vp (\|x\|_1)} \}$ a
$(\vp ,\PP)$\emph{-divergent volume-jet} if
$\{x\}$ is a toric jet for a unit jet $\{w\}$ framed by $\bar w$
that is adapted to $(\S,\C,\R(\PP))$, such that
{$\lim_{j\to \infty}\|x\|_1 = \infty\,$}.
\label{d:vjss}}
\noindent
As we show next, using this framework
further yields the estimate \eref{e:lambdaw}
(which, as outlined in Section~\ref{s:prop2}, is
the first step in proving {\eref{dfn:avx}).}
\lmm{\label{r:noindicatoratboundary}
{Setting $\vp (\rho) = e^{\rho}$, there exists $\delta>0$
such that for any {frame} $\bar w$, $r \in \R(\PP)_{\bar w^{(1)}}$
and
$(\vp ,\PP)$}-divergent volume-jet $(v,x)$ framed by $\bar w$, eventually,
\begin{equ}
\lambda_r(x) \ge \Lambda_r^{(\vvi)}(\xvi) \ge \delta \lambda_r(x)\,.\label{e:divlambdasuper}
\end{equ}}

\begin{proof} {Letting $\xi(j):=j! j^{-j}$ for $j \in \mathbb N$ and $\xi(0)=1$,
we set
\[
\delta_r := \prod_{i=1}^d \xi ( (c^r_\inn)_i) > 0 \;.
\]
As mentioned before, comparing \eref{e:jrates} and \eref{e:asrates}
one gets the first inequality of \eref{e:divlambdasuper} for any
$x \in (v^{-1}\mathbb N_0)^d$, $v \ge 1$. Further, the
ratio $\Lambda_r^{(v)}(x)/\lambda_r(x)$ is non-decreasing
in each $v x_i$ and equals $\delta_r$ when $v x = c_\inn^r$.
Thus, setting $\delta=\min_r \delta_r$ it suffices to show that
for $r \in \R(\PP)_{\bar w^{(1)}}$ and a
$(\vp ,\PP)$-divergent volume-jet {$\{( v, x)\}$} framed by $\bar w$,
we eventually have $v x_i \ge (c_\inn^r)_i$. This
trivially holds if $(c_\inn^r)_i = 0$,
so the proof is complete upon showing that, {along $\{(v,  x)\}$\,,}
\begin{equation}\label{eq:left-to-do}
i \in \supp\{c_\inn^r\} \quad \Longrightarrow \quad
\lim_{j \to \infty} \{\log v + w_i \log \theta\} = \infty \,.
\end{equation}
Since $(\log \|\pi_\PP x\|_1)/(\log \theta) \to \max_i \{\bar w^{(1)}_i\} =: \psi$
and
both $\|x\|_1$ and $\theta$ diverge,
we have $\psi \ge 0$. Further, $w_i \to \bar w^{(1)}_i$ is finite
and $\log v \ge \log \vp (\|x\|_1) =  \|x\|_1$ so
\eref{eq:left-to-do} clearly holds whenever $\psi>0$. In case
$\psi=0$ the vector $\bar w^{(1)} \in S^\PP$ has non-positive
coordinates, so $\bar w^{(1)}_{i'} \le -1/\sqrt{d}$ for some $i' \in \PP$.
Since {$\lim_j w = \bar w^{(1)}$}, it then follows that
eventually $w_{i'} \le - 1/\sqrt{2d} =: -\zeta$. Since $i' \in \PP$
and $v x \in {\mathbb N}^d(\PP)$ this implies that
$v \ge \theta^{-w_{i'}} \ge \theta^{\zeta}$. Recall
\rref{r:ired-in} that some $r' \in \R(\PP)$ has $c_\inn^{r'} = 0$,
hence $\langle \bar w^{(1)}, \pi_\PP c_\inn^r \rangle \geq 0$
for any $r \in \R(\PP)_{{\bar w}^{(1)}}$. That is, when
$i \in \supp\{c_\inn^r\}$ we have that $w_i \to \bar w^{(1)}_i = 0$
and as $\log v \ge \zeta \log \theta$, we recover \eqref{eq:left-to-do}
and with it, complete the proof.}
\end{proof}

Finally, adapting \cite[Defn. 6.8, 6.15]{gopal13}, each possible frame
within $S^\PP$, induces two key indices (classifications) for
reactions $r \in \R(\PP)$.
\dfn{
For $\PP \subseteq \S$ and ordered \abbr{ons} $\bar w \subset S^\PP$:
\newline
(a) Let $\super_1 = \R(\PP)_{\bar w^{(1)}}$
and
define
for $k \ge 2$
the nested subsets $\super_{k}$ of reactions
having $c_\inn^r$ in the $\bar w^{(k)}$-maximal
subset of $\{ \pi_\PP c_\inn^r : r \in \super_{k-1} \}$\,.
\newline
(b) The \emph{level} $\ell$ within $\bar w$ of
$r \in \R(\PP)$ having $\supp\{ c^r_{\rm out} \} \subseteq \PP$
is $\ell := \inf \{ k : \langle \bar w^{(k)} , {\pi_\PP} c^r\rangle \neq 0 \}$
(with $\ell=\infty$ when no such $k$ exists),
setting $\ell = 1$ if $r$ has some product species outside $\PP$.
\label{dfn:super-level}
}

\subsection{The dominance of dissipative reactions}
\label{s:tj}

Turning to the proof of \eref{dfn:Prrp}, we first bound
(in the setting of \lref{r:noindicatoratboundary}),
the contribution to the Lyapunov term {when
$r' \in \R(\PP)_{\bar w^{(1)}}$ and $\supp\{c_\outt^{r'}\} \not \subseteq \PP$,
allowing us thereafter to simultaneously treat
such reactions and those in $\R(\PP)_{\bar w^{(1)}}$
with $\supp \{c_\outt^{r'}\} \subseteq \PP$, ultimately using their
negative contribution to dominate any positive term in $a^{(v)}(x)$
from \eref{dfn:avx}.}
\lmm{For $\vp (\rho) = e^{\rho}$ and any
$(\vp ,\PP)$-divergent volume-jet $(v,x)$ framed by $\bar w$:
\newline
(a). If $r \in \R(\PP)_{\bar w^{(1)}}$ and
$\supp\{c_\outt^r\} \not \subseteq \PP$, then
\begin{equation}\label{e:e-l/22}
\limsup_{j \to \infty} \Big\{ \frac{h_r^{(v)}(x)}{\log \theta} \Big\} <  0 \,.
\end{equation}
(b). If $r \in \R(\PP)$ has $\supp \{c_\outt^r\} \subseteq \PP$ then
$\kappa_r = \infty$ iff\/
$r$ has finite level $\ell$
and
$\lim_j \theta{{{}}}^{\beta^{(\ell)}{{{}}}} = \infty$\,.
\label{l:e-1/22}
}
\begin{proof} (a). Recall that $\h{v}{x} = \langle \nablavU, c^r \rangle$, so
setting $\alpha_r := \langle \Ind {\PP^c} , c_\outt^r \rangle > 0$
and $\eta_r := \langle \Ind {\PP^c}  \log c_\outt^r , c_\outt^r \rangle$
which is finite, we have from \eqref{dfn:gradU-P} that
\begin{equation}\label{e:hrv-theta}
\frac{h_r^{(v)}(x)}{\log \theta} = \frac{\eta_r}{\log \theta} +
\langle w, \pi_\PP c^r \rangle  - \alpha_r \frac{\log v}{\log \theta}  \,.
\end{equation}
Because $\theta=\theta(j) \to \infty$,
the first term on the \abbr{rhs} decays to zero and the second term
converges to $\langle \bar w^{(1)} , \pi_\PP c^r \rangle$.
While proving \eqref{eq:left-to-do} we have seen that if
$\psi := \max_i \{ \bar w^{(1)}_i \} > 0$, then
$\log v \ge \|x\|_1$ (for the $\vp $-divergent volume-jet),
results with $(\log v)/(\log \theta) \to \infty$ and
consequently \eqref{e:e-l/22} holds. In case $\psi=0$ we have
shown in that same proof that
$(\log v)/(\log \theta) \ge \zeta > 0$ along the
divergent volume-jet and further that
$\langle \bar w^{(1)} , \pi_\PP c_\inn^r \rangle = 0$
when $r \in \R(\PP)_{\bar w^{(1)}}$. Recall
that $c^r = c_\outt^r - c_\inn^r$ with $c_\outt^r \in \RR^d$
and $\psi=0$ amounts to $-\bar w^{(1)} \in \RR^{d_\PP}$.
Thus, in this setting $\langle  \bar w^{(1)} , \pi_\PP c^r \rangle \le 0$,
which by \eqref{e:hrv-theta} recovers \eqref{e:e-l/22}.
\newline
(b). If $r \in \R(\PP)$ has $\supp\{c_\outt^r\} \subseteq \PP$ then
$\h{v}{\cdot}=h_r(\cdot)$ is independent of $v$ and in
\eqref{e:hrv-theta} we have $\alpha_r=\eta_r=0$. Further,
recall \dref{dfn:toric-jet} that $\{w\} \subset {\rm Co}(\bar w)$,
so if $r$ has infinite level
then $h_r(\cdot) = 0$,
while if it has finite level $\ell$ within $\bar w$,
then by \eqref{dfn:u-jet}, along the divergent volume-jet
\begin{equ}
\lim_{j \to \infty}
\frac{h_r^{(v)}(\xvi)}{\beta^{(\ell)}\log {\theta{{{}}}}}
= \langle \bar w^{(\ell)}, \pi_\PP c^r \rangle \ne 0 \,,
\label{e:logbetai}
\end{equ}
from which the stated criterion for divergence of
$|h_r^{(v)}(x)|$ follows.
\end{proof}

{
Our next result shows that $L_{r'}^{(v)}(x) \to -\infty$
for any dissipative {$r' \in \R(\PP)_{\bar w^{(1)}}$} with $\kappa_{r'} = \infty$
(see Sec. \ref{s:prop2} for explanation about the Lyapunov domination).
}

\lmm{\label{l:innerproductdomination}For $\vp (\rho)=e^\rho$, if $r \in \R(\PP)_{\bar w^{(1)}}$
with $\kappa_r = \infty$ is dissipative
for a $(\vp ,\PP)$-divergent volume jet $(v,x)$ framed by $\bar w$, then
\begin{equ}\Label{e:linnerproductdomination}
\lim_{j \rightarrow \infty} \UQv1 = -\infty~.
\end{equ}
}
\begin{proof} By \dref{d:vjss} the toric jet $\{x\}$ is adapted
to $(\S,\C,\R(\PP))$. Hence, if $\supp\{c_\outt^r\} \subseteq \PP$
and $r \in \R(\PP)$ is dissipative, then by \dref{d:dissipative}
and \eqref{e:grad-U},
\begin{equ}
h_r^{(v)}(x) = h_r(x)
= (\log \theta) \langle w, \pi_\PP c^r \rangle < 0 \,, \qquad \forall j \,.
\end{equ}
Since $\kappa_r=\infty$ it follows that $h^{(v)}_r(x) \to -\infty$
as $j \to \infty$, which by part (a) of \lref{l:e-1/22} applies
also when $\supp \{c_\outt^r\} \not \subseteq \PP$. Fixing $\gamma<\infty$,
since $\epsilon(x)$ of \lref{l:taylor} is uniformly bounded, we thus
have that for all $j$ large enough,
\[
L^{(v)}_r (x) = U(x) (Q^{(v)}_r(x) -1 )
\le U(x) \Big[ e^{-\frac{\gamma}{U(x)}} -1 \Big]
\le  - \gamma + \frac{\gamma^2}{2 U(x)} ~,
\]
(as $e^{-y} \le 1 - y + \frac{y^2}{2}$ for $y \in \RR$).
Recalling from \dref{d:vjss} that $\|x(j)\|_1 \to \infty$ and
consequently $U(x(j)) \to \infty$, we complete the
proof by
taking
$j \to \infty$ followed by $\gamma \to \infty$.
\end{proof}

We plan to show that if $r \in \R(\PP)$ has
$L_r^{(v)}(x) > 0$ along some $(\vp ,\PP)$-divergent
volume-jet $\{(v,x)\}$ for $\vp (\rho)>e^\rho$, then
\eref{dfn:Prrp} holds for a $\{w\}$-dissipative
$r' \in \R(\PP)_{\bar w^{(1)}}$.
To this end,
we first introduce the \abbr{crn} $\C_{\bar w^{(1)},\PP}$
in which necessarily $\supp \{c_\outt^r \} \subseteq \PP$ (or else
by \lref{l:e-1/22} (a) and \lref{l:innerproductdomination}
eventually $L_r^{(v)}(x) < 0$).
\dfn{For
$\PP \subseteq \S$ and $u \in S^\PP$ let
$(\S,\C_{u,\PP},\R(\PP))$ denote the \abbr{crn} obtained
upon restricting $c_\outt^r$ to $\RR^d(\PP)$
for any
$r \notin \R(\PP)_u$.
\label{dfn:mod-crn}}
\begin{remark}
Of course $\C_{u,\PP}=\C$ when $\PP=\S$. More generally, this
modification never affects $\{c_\inn^r\}$, hence neither the
rates $\Lambda^{(v)}_r(\cdot)$ nor the sets $\{ \R(\PP)_w,
w \in S^{\PP}\}$, or $\super_k$ of \dref{dfn:super-level}.
Further,
the \abbr{crn}
$(\S,\C_{u,\PP},\R(\PP))$ remains $u$-strongly
$\PP$-endotactic (see \dref{d:Pendo}) and thus also
$w(j)$-strongly $\PP$-endotactic, for $j$ large enough
and any unit jet $\{w(j)\}$ whose frame starts at
$\bar w^{(1)}=u$.
\label{rmk:ext-gopal}
\end{remark}
Comparing our \dref{dfn:toric-jet} and \dref{dfn:super-level} with the
corresponding definitions of \cite[Sec. 6]{gopal13}, it is easy to verify
that \cite[Thm.~6.11, Lemmas 6.7, 6.10, 6.19]{gopal13} apply in our setting
as does \cite[Prop.~6.20.1]{gopal13} (for draining reactions),
even for the modified \abbr{crn} of \dref{dfn:mod-crn}. We next adapt
to the latter setting, those conclusions of
\cite[Lemma 6.22, Prop.~6.24]{gopal13} that we shall use in the sequel.

{
\lmm{Fix a strongly endotactic \abbr{crn} \net. Consider the corresponding \abbr{crn} of \dref{dfn:mod-crn} and ordered
\abbr{ons} $\bar w$ of length $\ell'$ starting at $\bar w^{(1)}=u$. Then
\newline
(a). If $\supp \{ c_\outt^{r} \} \subseteq \PP$,
$\langle \bar w^{(k)}, \pi_\PP c^{r} \rangle=0$ for
$k < \ell'$ and
$\langle \bar w^{(\ell')}, \pi_\PP c^{r} \rangle > 0$,
then $r \notin \super_{\ell'}$\,.
\newline
(b). Some $r' \in \super_{\ell'}$ has
$\supp\,\{ c_\outt^{r'}\} \not \subseteq \PP$ or
$k \mapsto \langle \bar w^{(k)}, \pi_\PP c^{r'} \rangle$
not identically zero, with a negative first non-zero term.
\label{l:gopal-lem}}
\begin{proof}
Since $k \mapsto \super_{k}$ are nested sets, it suffices to rule out
that respectively:\\
(a$^{\prime}$). Some $r \in \super_{\ell'}$ has
$\supp\{c_\outt^{r}\} \subseteq \PP$, $\langle \bar w^{(k)}, \pi_\PP c^{r} \rangle = 0$ for $k < \ell'$ and $\langle \bar w^{(\ell')}, \pi_\PP c^{r} \rangle > 0$.\\
(b$^{\prime}$).
Each $r \in \super_{\ell'}$ has $\supp\{c_\outt^{r} \}\subseteq \PP$
and $\langle \bar w^{(k)}, \pi_\PP c^{r} \rangle = 0$ for all $k \le \ell'$.
\newline
Further, the modification of \dref{dfn:mod-crn} neither
affects $\super_{\ell'}$ nor the value of $c^r$ for reactions in
$\super_{\ell'} \subseteq \R(\PP)_{u}$ (see \rref{rmk:ext-gopal}),
so it suffices to rule out (a$^{\prime}$) and (b$^{\prime}$) for
$(\S,\C,\R(\PP))$ and the given \abbr{ons} $\bar w$.
To this end, consider a unit jet
$\{w\}$ framed by $\bar w$, adapted to
$(\S,\C,\R(\PP))$ and having $\beta^{(\ell')}(j) > 0$ for all $j$.
Recall \cite[Thm.~6.11]{gopal13} that
$\R(\PP)_{w(j)} = \super_{\ell'}$ eventually in $j$. Thus, by
\eref{dfn:u-jet}, our assumptions (a$^{\prime}$) resp.~(b$^{\prime}$)
imply that for all large enough $j$, respectively:\\
(a$^{\dagger}$). There exists a $w(j)$-explosive $r \in \R(\PP)_{w(j)}$ of level $\ell'$.\\
(b$^{\dagger}$). The collection $\R(\PP)_{w(j)}$ consists of only $w(j)$-null reactions.
\newline
To conclude, note that (a$^{\dagger}$) and (b$^{\dagger}$) contradict having a
strongly $\PP$-endotactic $(\S,\C,\R(\PP))$.
\end{proof}}

Similarly to
\cite[Prop.~6.26]{gopal13}, we proceed via a pair of lemmas that
establish \eref{dfn:Prrp} for $(\S,\C_{\bar w^{(1)},\PP},\R(\PP))$
by bounding from
below the asymptotic behavior of the Lyapunov and monomial terms,
as in cases (a) and (b) at the end of Sect. \ref{s:prop2}, that
correspond to $\kappa_r<\infty$ and $\kappa_r =\infty$, respectively.

\lmm{[Lyapunov domination] For $\vp (\rho) = e^\rho$ and
the \abbr{ons} $\bar w$ for $\PP \subseteq \S$, consider
the \abbr{crn} $(\S,\C_{\bar w^{(1)},\PP},\R(\PP))$ and
a $(\vp ,\PP)$-divergent volume jet $(v,x)$ for it, framed by $\bar w$.
Then, for any $r \in \R(\PP)$ with $\supp \{c_\outt^r\} \subseteq \PP$
and $\kappa_r < \infty$, the domination \eqref{dfn:Prrp} holds for
some dissipative $r' \in \R(\PP)_{\bar w^{(1)}}$.
\label{l:innerproductdomination3}}

\begin{proof} Let $\ell$ denote the level of $r \in \R(\PP)$
within the frame $\bar w$, if finite, whereas if
the level of $r$ is infinite, set
$\ell = d_\star+1$ and $\beta^{(\ell)} \equiv 0$.
Since $\supp\{c_\outt^r\} \subseteq \PP$ we have from
\lref{l:e-1/22} (b) that
$\lim_j \theta^{\beta^{(\ell)}} < \infty$\,.
For any divergent volume-jet $\beta^{(1)} \to 1$,
hence $\lim_{j} \theta{{{}}}^{\beta^{(1)}{{{}}}} = \infty$,
$\ell \ge 2$ and in view of \eqref{dfn:u-jet}
there exists $1 \le \ell' < \ell$ such that
\begin{equ}
\lim_{j \rightarrow \infty} \theta{{{}}}^{\beta^{(\ell')}{{{}}}} = \infty \anddd \lim_{j \rightarrow \infty} \theta{{{}}}^{\beta^{(\ell'+1)}{{{}}}} < \infty\,.
\label{e:ellprime}
\end{equ}

For the sub-frame $\{\bar w^{(1)}, \dots, \bar w^{(\ell')}\}$,
\lref{l:gopal-lem}(b)
yields $r' \in \super_{\ell'} \subseteq
\R(\PP)_{\bar w^{(1)}}$
of level $\ell_\star \le \ell'$ such that
either $\supp\{c^{r'}_\outt\} \not\subseteq \PP$ or
$\langle \bar w^{(\ell_\star)}, \pi_\PP c^{r'} \rangle < 0$.
Since
$\lim_j \beta^{(k)}/\beta^{(k+1)} = \infty$ for any $k \ge 1$,
such $r'$ must also be $\{w\}$-dissipative.
Proceeding to establish \eqref{dfn:Prrp},
by \lref{l:taylor} combined with $e^{|x|}-1 \leq 2|x|$ for $|x| < 1$
and $\h{v}{x}/U(x) \to 0$, we have for $j$ large enough
$|\UQv1| \leq 2 (|\epsilon(x)| + |{\h{v}{\xvi}}|) \leq C^{-1}$, hence
\begin{equ}
{\prrr}
\geq C \tg{\pi_\PP(c^{r'}_{\inn}-c^{r}_{\inn})} |L_{r'}^{(v{{{}}})}(\xvi)| ~.
\Label{e:pinnerproductdomination3}
\end{equ}
As $r' \in \super_{\ell'}$ and $c_\inn^r \in \C_\inn(\PP)$ we have
from \cite[Lemma 6.10.2]{gopal13}
that for any $k_\star \le \ell'+1$, $\delta>0$ and all $j$ large enough
\begin{align}
\langle w(j), \pi_\PP (c_\inn^{r'} - c_\inn^r) \rangle \ge &
\sum_{k = k_\star}^{d_\star} \beta^{(k)}(j) \langle \bar w^{(k)} ,
\pi_\PP (c_\inn^{r'} - c_\inn^r)
\rangle \nonumber \\
&\ge \beta^{(k_\star)} (j)
[\langle \bar w^{(k_\star)} , \pi_\PP (c_\inn^{r'} -c_\inn^r) \rangle - \delta]
\,.
\label{eq:new-6-10-2}
\end{align}
Taking $k_\star=\ell'+1$ (where if $\ell'=d_\star$ then
$\langle \bar w^{(k)}, \pi_{\PP}( c_\inn^{r'} - c_\inn^r ) \rangle \ge 0$
for all $k \le d_\star$ hence the \abbr{lhs} of \eref{eq:new-6-10-2} is
non-negative), we deduce from \eref{e:ellprime} that
$\tg{\pi_\PP(c^{r'}_{\inn}-c^{r}_{\inn}) }$ is uniformly (in $j$) bounded below.
The proof
is thus complete upon showing
that $\kappa_{r'}=\infty$, as then
$|L_{r'}^{(v)}(x)| \to \infty$ by \lref{l:innerproductdomination}.
Indeed, since $r' \in \R(\PP)_{\bar w^{(1)}}$
from \lref{l:e-1/22}(a) we have that $\kappa_{r'}=\infty$
if $\supp\{c^{r'}_\outt \} \not \subseteq \PP$, whereas
if $\supp\{c^{r'}_\outt\} \subseteq \PP$ then
$r'$ of finite level $\ell_\star \le \ell'$ has $\kappa_{r'}=\infty$
in view of the \abbr{lhs} of \eref{e:ellprime} and
\lref{l:e-1/22}(b).
\end{proof}

\lmm{[Monomial domination]
For $\vp (\rho) = e^\rho$ and
\abbr{ons} $\bar w$ for $\PP \subseteq \S$, consider
the \abbr{crn} $(\S,\C_{\bar w^{(1)},\PP},\R(\PP))$ and
a $(\vp ,\PP)$-divergent volume jet $(v,x)$ for it, framed by $\bar w$.
Then, for any $\{w\}$-explosive
$r \in \R(\PP)$ with $\supp \{c_\outt^r\} \subseteq \PP$
and $\kappa_r = \infty$, the domination \eqref{dfn:Prrp} holds for
some dissipative $r' \in \R(\PP)_{\bar w^{(1)}}$.
\label{l:monomialdomination}}
\begin{proof} By \lref{l:e-1/22}(b) here $r$ has finite level $\ell'$
within $\bar w$ for which the \abbr{lhs} of \eref{e:ellprime} holds.
Further, with $\{w(j)\}$ adapted to $(\S,\C_{\bar w^{(1)},\PP},\R(\PP))$ we
deduce from \cite[Prop.~6.20.1]{gopal13} that since
$\langle w(j), \pi_\PP c^r \rangle$ is positive for $j$ large,
$\langle \bar w^{(\ell')}, \pi_\PP c^r \rangle$ must also be positive,
hence \lref{l:gopal-lem}(a) yields that
$r \notin \super_{\ell'}$.
Recall the proof of \lref{l:innerproductdomination3}, that
there exists $\{w\}$-dissipative
$r' \in \super_{\ell'} \subseteq \R(\PP)_{\bar w^{(1)}}$.
{In particular,
$\langle \bar w^{(\ell')}, \pi_\PP (c_\inn^{r'} - c_\inn^r) \rangle$
is positive, so considering \eqref{eq:new-6-10-2} for $k_\star=\ell'$
and small $\delta>0$, for $j$ large enough we bound the monomial
term of \eqref{dfn:Prr} by
\begin{equs}
\tg{\pi_\PP (c^{r'}_\inn - c^r_\inn)} \ge
\Big(\theta^{\beta^{(\ell')}}\Big)^{\delta} \,.
\label{e:pmonomialdomination1}
\end{equs}
Further, the $\{w\}$-dissipative
$r' \in \R(\PP)_{\bar w^{(1)}}$ has level $\ell_\star \le \ell'$ and $\kappa_{r'}=\infty$, hence
\begin{equs}
K_{r'} := \lim_{j \to \infty}
\frac{h_{r'}^{(v)}(\xvi)}{\beta^{(\ell')}\log {\theta{{{}}}}}
\label{e:hrp-div-rate}
\end{equs}
is strictly negative
(see \eqref{e:e-l/22} for $\supp\{c_\outt^{r'}\} \not \subseteq \PP$
and \eqref{e:logbetai} otherwise).
The $\{w\}$-explosive $r$ has
level $\ell'$ and $\kappa_r=\infty$ hence by \eref{e:logbetai} it
satisfies \eref{e:hrp-div-rate} for some $0<K_r<\infty$. Recall
\eref{e:ellprime} that $\beta^{(\ell')}\log \theta$ diverges along our
jet $\{w\}$. Hence, by Lemma \ref{l:taylor}, for any $s > K_r$
and $\gamma \in (0,1)$ such that $\gamma s < - K_{r'}$, the
corresponding Lyapunov term
is eventually bounded below
by
\begin{equs}
\frac{1-\qqq}{\qrq-1} \ge \frac{1-\Big(\theta^{-s \beta^{(\ell')}/U(x)}
\Big)^{\gamma} }{\theta^{s \beta^{(\ell')}/U(x)} -1} \ge
\gamma \theta^{-s \beta^{(\ell')}/U(x)}
\label{e:plyapunovdomination-new}
\end{equs}
(where the second inequality follows from
$1-\xi^\gamma \ge \gamma (1-\xi)$ which holds for any $\xi,\gamma \in (0,1)$).
With $U(x) \to \infty$, the \abbr{rhs} of \eref{e:pmonomialdomination1}
dominates the \abbr{rhs} of \eref{e:plyapunovdomination-new} and the
divergence of $P_{r,r'}^{(v)}(x)$ of \eref{dfn:Prr} follows.}
\end{proof}

\subsection{Proof of \pref{t:thm}}
\label{s:proofofprop}

By \eref{dfn:avx}, \pref{t:thm} will hold if we can find
$\rho<\infty$ such that for any $\rho' < \infty$,
\[
\vp (\rho') := \sup \{ v : \sup_{x \in \AN^v} \{a^{(v)}(x)\} > 0 \}
< \infty \,.
\]
Assume to the contrary, that there exist $\rho'_j > \rho_j \uparrow \infty$, $\vvij \to \infty$ as $k \to \infty$
and $\xvij \in \mathcal{A}_{\rho_j,\rho'_j}^{\vvij}$ such that
$a^{\vvij}(\xvij) > 0$ for all $j,k \in \mathbb{N}^2$.
Then, for any desired increasing $\vp (\cdot)$, upon choosing $k=k_j$
large enough, we extract a sequence $\{(v(j),x(j))\}$ such that
$v(j) x(j) \in {\mathbb N}_0^d$, $\|x(j)\|_1 \to \infty$ and
\begin{equ}
a^{v(j)} (x(j)) > 0  \,, \quad v(j) > \vp (\|x(j)\|_1) \,,
\qquad \forall j \in \mathbb{N}\,.
\label{e:absurdum}
\end{equ}
Since $d<\infty$, there must be some $\PP \subseteq \S$
such that $v(j) x(j) \in {\mathbb N}^d(\PP)$
along some infinite sub-sequence. Also, as
$\|x(j)\|_1 \to \infty$,
upon restriction to $(\S,\C,\R(\PP))$ we have that
$\theta(x(j)) \to \infty$ (see \eqref{dfn:topic-p}),
and our \dref{dfn:toric-jet} of unit jet and toric jet
then coincide with those of \cite{gopal13}. Hence,
by \cite[Lemma 6.7]{gopal13} we
extract a sub-sub-sequence $(v(j),x(j))$ satisfying all of the
above, for which in addition $\{x(j)\}$
is a toric jet for a unit jet $\{w(j)\}$
framed by some $\bar w$.
Finally, in view of \cite[Lemma 6.19]{gopal13},
there exists a further sub-sub-sub-sequence $\{x(j)\}$
which is adapted to $(\S,\C,\R(\PP))$ (note that
$\supp\{c^r_{\outt}\} \not \subseteq \PP$ has nothing to do with
the choice of $\{w(j)\}$). In conclusion, we have a
$(\vp ,\PP)$-divergent volume-jet $\{(v,x)\}$ satisfying
\eref{e:absurdum}, where we are free to choose $\vp (\rho)$ and
only $r \in \R(\PP)$ is to be considered in \eref{dfn:avx}.
Fixing $\{(v,x)\}$ and in particular its frame $\bar w$,
we may and can move to the \abbr{crn}
$(\S,\C_{\bar w^{(1)},\PP},\R(\PP))$ of \dref{dfn:mod-crn}.
Indeed, recall \rref{rmk:ext-gopal} that
this does not affect
the rates $\Lambda_r^{(v)}(x)$, while for $v \ge \sup_r \|c_\outt^r\|_\infty$ and
$x \in \RR^d(\PP)$ it may only increase $L_r^{(v)}(x)$,
by setting to zero some negative
contributions $(v^{-1} c^r)_i [\log (v^{-1} c^r)_i - 1]$ to
$U(x+v^{-1} c^r)$ from $i \in \supp \{ c^r_\outt \} \setminus \PP$.
As explained before, {in the new \abbr{crn}
$\supp \{c_\outt^r\} \not \subseteq \PP$
requires $r \in \R(\PP)_{\bar w^{(1)}}$ and further
sub-sampling our divergent volume-jet to make it adapted
to $(\S,\C_{\bar w^{(1)},\PP},\R(\PP))$, we proceed
as outlined in Section~\ref{s:prop2} to show that
on the latter \abbr{crn}, having \eref{e:absurdum}
leads to a contradiction. Indeed, consider $r \in \R(\PP)$,
whose contribution to \eref{e:absurdum} is eventually positive
(for the modified reactions of $\C_{\bar w^{(1)},\PP}$).
That is, having $\UQv1>0$ for all $j$ large.
By \lref{l:taylor} this requires $\h{v}{x} + \epsilon(x)>0$, which in
view of \lref{l:e-1/22} (a) implies
that $\supp \{c_\outt^r \} \subseteq \PP$.
With $\{x\}$ adapted, this yields, as in
the proof of \lref{l:e-1/22} (b), that
$|h_r^{(v)}(x)| \to \kappa_r$ when $j \to \infty$ (see
\eref{e:logbetai}), and further that $\kappa_r = \infty$
is possible only for a $\{w\}$-explosive reaction.
For both $\kappa_r < \infty$ and $\kappa_r=\infty$ we now have
\eref{dfn:Prrp} for some dissipative $r' \in \R(\PP)_{\bar w^{(1)}}$
(see \lref{l:innerproductdomination3} and \lref{l:monomialdomination},
respectively). As \eref{dfn:Prr} is a consequence of
\lref{r:noindicatoratboundary}, it follows that $a^{(v)}(x) \leq 0$ along
$\{(v,x)\}$, in contradiction with \eref{e:absurdum}.
}

\section{Proof of \twf}
\label{s:wf}

\twf~ is proved in \cite[§ 6]{fw98} for a uniformly elliptic
diffusion on a compact $d$-dimensional manifold, when the
driving Brownian motion has been scaled by $\epsilon$. Recall that
such a diffusion satisfies an \abbr{LDP} with rate $v:=\epsilon^{-2}$
and its good rate function
is zero iff $x'(t) = b(x(t))$ starting at $x(0)=x_0$. We have
here the analogous \abbr{ldp} of \tref{t:largedev}, whose good rate
function is zero iff $z(t)$ solves the \abbr{ode} \eref{e:ma}
(see \rref{r:zero-rate}). Further, with our Assumptions
\ref{a:3} and \ref{a:2} replacing \cite[Condition A, § 6.2]{fw98}
and \cite[§ 6.5]{fw98}, respectively, we merely adapt the proof
in \cite[§ 6]{fw98}, where the stated results are established
from \cite[Lemmas 6.1.1--6.1.9]{fw98}. Specifically,
for \eref{e:twf1} and \eref{e:twf2} which concern only
the dynamics of $t \mapsto X_t^v$ within the compact $\mathcal D$,
it suffices that we prove the weaker version \lref{l:11fw} of
\cite[Lemma 6.1.1]{fw98} within $\mathcal D$, {and the modification \lref{l:14fw}
of \cite[Lemma 6.1.4]{fw98},}
while tackling the degeneracy of
$\{X_t^v\}$ on $\partial \RR^d$.
Indeed, \lref{l:11fw} {and \lref{l:14fw}} suffice for
establishing \cite[Lemmas 6.1.2 and 6.1.4]{fw98} respectively. Furthermore, the local
Lipschitz continuity of the quasi-potential
is never used in the proof of \eref{e:twf1} and \eref{e:twf2},
while \cite[Lemma 6.1.3]{fw98}
can be bypassed (since it is only used for proving \cite[Lemma 6.1.4]{fw98}).
The \abbr{LDP} and \cite[Lemmas 6.1.1-6.1.4]{fw98},
together imply \cite[Lemmas 6.1.5-6.1.9]{fw98}, containing the fundamental transition times estimates for the establishment of \cite[Lemmas 6.2.1, 6.2.2]{fw98}, proving that $\mathcal V_{\mathcal D}$ is the relevant functional
for the estimation of transition probabilities between $K_i$'s. The
combination of these results finally yields \eref{e:twf1} and \eref{e:twf2}
as explained in the proofs of \cite[Thms.~6.5.1,~6.5.3]{fw98}.
We thus proceed to state and prove the {adaptations of} \cite[Lemmas 6.1.1 and 6.1.4]{fw98}
to the current setting.

\lmm{
For $\mathcal D \subset \RR^d$
as in \twf~
there exist $\kappa \ge 1$, $\epsilon>0$ and $C(t) \to 0$
(as $t\to 0$), such that for any
$x,y \in {\mathcal D}$ with $\|x-y\|_1 < \epsilon$,
there exists a path $z(\cdot) \subset {\mathcal D}$,
of length $t = \kappa \|x-y\|_1$ with
$I_{x,t}(z) \le C(t)$ and $z(t)=y$.
\label{l:11fw}
}
\begin{proof} By the continuity of $\lambda_r(\cdot)$ of \eqref{e:asrates}
on $\mathcal D$ compact,
$\bar \lambda :=
\max_{r \in \R, x \in \mathcal D} \{\lambda_r(x)\}$ is finite.
Further, since ${\rm Co} \{c^r\}_{r \in \R} = \Rr^d$ the sets $Q_\R(\xi)$ are non-empty and
$$
 \bar q := {e} \vee
 \max_{\|\xi\|_1 \le 1} \min \{ \|q\|_\infty : q \in Q_\R(\xi) \}
 < \infty \,.
$$
Setting
{${\bar c_\star} := {\sup_{r \in \R}} \{\|c^r_\inn\|_1 \}$
and $\gamma := \bar \lambda - \bar q + \bar q \log
(\bar q/\min_{r \in \R} \{k_r \wedge 1\})$ for the reaction constants
$k_r$ of \eqref{e:asrates},}
we then have for
any $z \in \mathcal D$ and $\|\xi\|_1 \le 1$ the bound
$$
L(\la(z),\xi) \le m \Big[ \gamma + \bar q {\bar c_\star} \big( \log \min_{i=1}^d \{z_i\} \big)_- \Big]
$$
on the Lagrangian of \eqref{e:lagrangian}. {Thus,
if $z \in AC_{0,t}\pc{\mathcal D}$ with $z(0)=x$
is such that $\|z'(s)\|_1 \le 1$ and
$\min_{i} \{z_i(s)\} \ge \beta s$, then for the rate function
of \eqref{dfn:Ixt},}
\begin{equs}
I_{x,t}(z) \leq  c(t) := {m \int_0^t \big[ \gamma + \bar q {{\bar c_\star}}
 (\log \beta s)_-  \big] \,\d s} \,.
\label{e:boundI}
\end{equs}
Similarly to \cite[Lemma 2.1]{SW2}, \aref{a:2} implies that
for some $\beta \in (0,1)$, $\epsilon \in (0,1/3)$
and $v^{(j)} \in \Rr^d$ with $\|v^{(j)}\|_1 \le 1$,
there exists a finite covering of $\mathcal D$ by
balls $\{\mathcal B_j\}$ such that
\begin{equs}
\min_{\widetilde x \notin \mathcal D}
\|x+s v^{(j)} - \widetilde x\|_\infty
\ge \beta s \,, \qquad \forall x \in \mathcal D \cap \mathcal B_j^\epsilon\,,
\;\;\; s \le \epsilon/\beta\,.
\label{e:boundDc}
\end{equs}
Fixing such a covering we set $\kappa = 1 + 2/\beta$. Suppose now
that $x \in \mathcal D \cap \mathcal B_j$
and $\|y-x\|_1 = \delta < \epsilon$ for some $y \in \mathcal D$.
Taking $t=t_1+t_2+t_3$ for $t_1=t_3=2\delta/\beta$ and $t_2=\delta$,
consider the continuous path from $x^{(1)}:=x$ to $x^{(4)}:=y$, composed of the
line segments between $x^{(1)}$, $x^{(2)}=x^{(1)}+t_1 v^{(j)}$,
$x^{(3)} = x^{(4)} + t_3 v^{(j)}$
and $x^{(4)}$. That is, $z^{(1)}(s)=x^{(1)} + s v^{(j)}$ for $s \in [0,t_1]$,
then $z^{(2)}(s)=x^{(2)} + \frac{s}{\delta} (y-x)$ for $s \in [0,t_2]$,
and finally, in reverse
$z^{(3)} (s) = x^{(4)} + s v^{(j)}$ for $s \in [0,t_3]$.
Since
$y \in \mathcal D \cap \mathcal B_j^\delta$ and $\delta \le \epsilon$,
it follows from \eqref{e:boundDc} that
$\min_i \{ z_i^{(\ell)}(s) \} \ge \beta s$ and
$z^{(\ell)}(s) \in \mathcal D$ for $\ell=1,3$ and $s \in [0,\delta/\beta]$.
The end points $x^{(2)}$ and $x^{(3)}$
of $z^{(2)}(\cdot)$, are $\delta$ apart and by the preceding, of
at least $2 \delta$ sup-distance from $\mathcal D^c$.
Consequently, {$\inf_{\xi \in \mathcal D^c}\|z^{(2)}(s)-\xi\|_1 \geq \delta$} and
$\min_i \{ z_i^{(2)}(s) \} \ge \delta \ge \beta s$ for $s \in [0,\delta]$.
By construction
$\|z'^{(\ell)}(s)\|_1 \le 1$ for $\ell=1,2,3$ and all $s$,
so in view of \eqref{e:boundI}
\[
I_{x,t}(z) = \sum_{\ell=1}^3 I_{x^{(\ell)},t_\ell}(z^{(\ell)})
\le \sum_{\ell=1}^3 c(t_\ell) =: C(t) \,,
\]
as claimed.
\end{proof}

\lmm{
Let
$\mathcal D_{- \delta}:= \mathcal D \setminus (\partial \mathcal D)^\delta$ with $\mathcal D$ as in \aref{a:2}. For some $C_\star(t) \to 0$, some
$\eta(\gamma,\kappa_\star,\mathcal D) > 0$,
any $\kappa_\star<\infty$, $\gamma > 0$ and $\delta \in (0,\eta)$,
if $T+I_{z_0,T}(z) \leq \kappa_\star$ for $z([0,T]) \subset \mathcal D$,
then there exists $\widetilde T \le T+ 3 \kappa \gamma$ and
$\tilde z([0,\widetilde T]) \subset \mathcal D_{- \delta}$ such that
$I_{\tilde z_0,\widetilde T}(\widetilde z) \leq I_{z_0,T}(z) + C_\star(\gamma)$ and
$\|\tilde z(0)-z(0)\|_1 + \|\tilde z(\widetilde T)-z(T)\|_1 \le 2 \delta$.
The same holds for $\mathcal D_{+\delta} := \mathcal D^\delta \cap \RR^d$ and $\mathcal D$,
instead of $\mathcal D$ and $\mathcal D_{-\delta}$, respectively.
\label{l:14fw}}

\begin{proof} From \cite[Lemma 2.1]{SW2} and \aref{a:2} we have \cite[Assmp. 2.1]{SW2} holding. Further, with $\bar \lambda$ finite, the path $z(\cdot)$ whose length and rate function are both bounded by $\kappa_\star$, makes at
most $J=J(\kappa_\star)$ transitions between the balls $\mathcal B_j$ in the covering of $\mathcal D$ (see \cite[Lemma 3.5]{SW2}).
Each of the monomials $\lambda_r(\cdot)$ of
\eqref{e:asrates}
is $c_{\mathcal D}$-Lipschitz continuous
on the compact $\mathcal D$ and
non-decreasing along any short path that
originates in a small enough neighborhood of {the set of zeroes of $\lambda_r(\cdot)$ in}
$\partial \RR^d$, and is directed inwards to $(\RR^d)^o$.
In particular, for some
$\nu>0$ and all $j$, \abbr{wlog}
the vectors $v^{(j)}$ in \eqref{e:boundDc} are such that
$\lambda_r(x+\alpha v^{(j)}) \ge \lambda_r(x)$
for any $\alpha \in [0,\nu]$ and $x \in \mathcal B_j$
for which $\lambda_r(x) \le \nu$.

Adapting \cite[Lemma 4.3]{SW2} we construct
for $\beta \in (0,1)$ as in the proof of \lref{l:11fw}
and some $\eta(\gamma,\kappa_\star,\mathcal D)>0$,
a path $\hat z \in (\mathcal D)_{-2\eta}$ with
$I_{{\hat z}_0, \hat T}(\hat z) \leq I_{z_0, T}(z) + 2\gamma$,
$\sup_t \|\hat z(t) - z(t)\|_1 \le \gamma$,
$\hat T \le T + \gamma$ and
$\|\hat z_0 - z_0\|_1 \le \eta' := 4 \eta/\beta$.
Specifically,
let $\hat z_0 = z_0 + \eta' v^{(i)}$ or
$\hat z_0 = z_0$ depending on whether $z_0 \in {\mathcal B}_i$
for $\mathcal B_i$ touching, or not touching, $\partial \mathcal D$, respectively.
Thereafter, $\hat z(\cdot)$ is parallel to $z(\cdot)$,
except that at the $k$-th time the path
$z(\cdot)$ transitions to a new ball $\mathcal B_j$ of the covering (that touches $\partial \mathcal D$),
a linear segment in direction $v^{(j)}$ is inserted
in $\hat z(\cdot)$ for duration $\eta_k = \eta' (3/\beta)^k$,
to keep it within ${\mathcal D}_{-2\eta}$.
With at most $J({\kappa_\star})$ transitions between different balls $\mathcal B_j$, taking $\eta>0$ small enough guarantees that the total contribution of time shifts to the length $\hat T$
of the path $\hat{z}$ be at most $\gamma$, and that
$\sup_s \|\hat{z}(s)-z(s)\|_1 \le \gamma$. Next, having
$I_{x,t}(x+s v^{(j)}) \le c(t)$, due to \eqref{e:boundI},
the rate contribution of all additional
linear segments is at most $\sum_k c(\eta_k) \le \gamma$
(for small enough $\eta>0$). Taking even smaller
$\eta>0$,
bounds by $\gamma$ (uniformly over all such path $z$),
the accumulated rate difference between pieces of
$\hat z(\cdot)$ and their parallels within $z(\cdot)$,
as soon as we show that for some $g_{\mathcal D} (\alpha) \to 0$
when $\alpha \to 0$,
\begin{equs}\label{e:mon-vj}
z(\cdot) \subset \mathcal B_j, \;\;  \alpha
\in [0, \nu/c_{\mathcal D}]
 \quad \Rightarrow \quad
I_{z_0,t}(z(\cdot)+\alpha v^{(j)}) \le I_{z_0,t}(z(\cdot)) + t g_{\mathcal D} (\alpha) \,.
\end{equs}
To this end, if $|\lambda_r -\hat \lambda_r| \le
c_{\mathcal D}  \alpha$ and
$\hat \lambda_r \ge \lambda_r$ whenever $\lambda_r \le \nu$, then
by \eqref{e:lagrangian}, for any $\xi \in \Rr^d$,
\begin{align*}
L(\hat \lambda,\xi) - L(\lambda,\xi)
\le
\|\hat \lambda - \lambda\|_1
+ \max_r \Big\{ \log \big(\frac{\hat \lambda_r}{\lambda_r}\big) \Big\}_-
\le m c_{\mathcal D} \alpha - \log (1- c_{\mathcal D}
\alpha/\nu) \,,
\end{align*}
hence denoting the \abbr{rhs} by $g_{\mathcal D}(\alpha)$
yields  \eqref{e:mon-vj} (see \eqref{dfn:Ixt}).

Now, fixing
$\delta \in (0,\eta)$, let $\tilde z(\cdot)$ be
$\hat z(\cdot)$ augmented by the initial/final
piece-wise linear path of \lref{l:11fw}, leading from
$\tilde z (0) := \arg \min_{z \in \mathcal D_{- \delta}}
\|z-z(0)\|_1$ to $\hat z(0)$ and from
$\hat z (\hat T)$ to $\tilde z(\tilde T) := \arg\min_{z \in \mathcal D_{- \delta}} \|z-z(T)\|_1$, respectively.
Since both $\|\tilde z(0)-\hat z(0)\|_1 \le \delta+\gamma$ and
$\|\hat z(\hat T) - \tilde z(\tilde T) \|_1 \le 2\gamma +
\delta$, taking $\eta \le \gamma \le \epsilon/3$ we have by
\lref{l:11fw} that the length of each augmented path is
at most $\kappa \gamma$ and its contribution to the
total rate does not exceed $C(3 \gamma)$. Finally, note that by
construction both end-points of these initial and final
pieces are in $\mathcal D_{-\delta}$, whereby the construction
of \lref{l:11fw} guarantees that their minimal distance
from $\partial \mathcal D$ be attained at one of their
end points, hence do not exceed $\delta$.
\end{proof}

While \eref{e:twf1} and \eref{e:twf2} involve only the process
$t \mapsto X^v_t$ within the compact $\mathcal D$, this is not
the case for \eref{e:twf4} which is established in
\cite[Thm.~6.6.2]{fw98} under the additional
assumption of a compact state space, which we lack here.
However, the latter proof applies for the stopping time
$\tau_{\pi,\rho} := \tau_\pi \wedge \sigma_\rho$ and the
non-random $C_\rho(\pi)$ obtained via \cite[Eqn. (6.6.1)-(6.6.2)]{fw98}
from $I^{(\rho)}_{x,t}(\cdot)$ of \eqref{dfn:Ixt} that
corresponds to $\lambda_r(x) \Ind{\widetilde K_\rho} (x)$,
with $\lambda_r(\cdot)$ of \eqref{e:asrates} and
$\widetilde K_\rho$ of \eqref{dfn:K-tilde}
(as the Markov jump processes $X_t^{v,\rho}$ from the proof of
\tref{t:largedev} are $\widetilde K_\rho$-valued and
satisfy the \abbr{ldp} with rate functions $I^{(\rho)}_{x,t}(\cdot)$).
For $\rho \ge \gamma$ and
$\cup_j K_j^\delta \subset \widetilde K_\gamma$
it is easy to verify that using
$I^{(\rho)}_{x,t}(\cdot)$ instead of $I_{x,t}(\cdot)$
amounts to replacing the quasi-potential $\VV(\cdot,\cdot)$
by $\VV_{\widetilde K_\rho}(\cdot,\cdot)$, with an
additional attractor of the dynamics at $(\widetilde K_\rho)^c$.
It is irrelevant that \aref{a:3} fails for this new attractor,
since it is outside $\pi$ hence the
transitions $(\widetilde K_\rho)^c \to K_j$ play
no role in $C_\rho(\pi)$.
By the same reasoning, the rate $I_{x,t}(z)$
of any path $z(\cdot)$ exiting $\widetilde K_\rho$
is part of the minimization yielding $C_\rho(\pi)$,
while those paths which are confined to $\widetilde K_\rho$
make exactly the same contribution to $C_\rho(\pi)$ and to $C(\pi)$.
Consequently, $C_\rho(\pi) \uparrow C_\infty(\pi) \le C(\pi)$ and
$v^{-1} \log \tau_{\pi,\rho} \to C_\infty (\pi)$ when $v \to \infty$
followed by $\rho \to \infty$. The compact sets $\widetilde K_\rho$
satisfy \aref{a:2}, so by \lref{l:11fw} the quasi-potential
$\VV(x,y)$ is everywhere finite. This implies that
$C(\pi)$ is finite, and thereby so is $C_\infty(\pi)$.
Considering \lref{t:exptight} for some $\beta > C_\infty(\pi)$
and $\rho \to \infty$, we thus conclude that
$v^{-1} \log \tau_\pi \to C_\infty(\pi)$, which translates
to \eref{e:twf4} provided $C_\infty(\pi) \ge C(\pi)$. The latter
is a direct consequence of our next lemma, showing that
$\VV(\widetilde K_\gamma,(\widetilde K_\rho)^c)
\to \infty$ as $\rho \to \infty$. Indeed, the second term on the
\abbr{RHS} of \cite[Eq.~(6.6.2)]{fw98} is independent of the addition of
$(\widetilde K_\rho)^c$ to the set of attractors (hence identical for $C(\pi)$ and $C_\rho(\pi)$), while every element over which the minimum is taken in
\cite[Eq.~(6.6.1)]{fw98} is either the same for $C(\pi)$ and
$C_\rho(\pi)$, or involves some transition $K_j \to (\widetilde K_\rho)^c$.
Since $\VV(\cdot,\cdot) \ge 0$, terms involving any such transition are
irrelevant when $\VV(\widetilde K_\gamma,(\widetilde K_\rho)^c) > C(\pi)$.

\lmm{Under \aref{a:1}, for any $\gamma$ finite,
\begin{equs}
\lim_{\rho \to \infty}  \inf_{t \ge 0} \{J_\gamma (t,\rho)\} = \infty \,, \quad
J_\gamma (t,\rho) := \inf_{\|x\|_1 \le \gamma} \;
\inf_{\{z(\cdot)~:~\sup_{s \le t} \|z(s)\|_1 > \rho\}} \, \{I_{x,t}(z)\} \,.
\label{e:quas-pot-inf}
\end{equs}
\Label{l:twf3}}
\begin{proof} The lower bound of the \abbr{ldp}
of \tref{t:largedev} for the open set
$\Gamma := \{ z : z(t) \in (\widetilde K_\rho)^c$ for some $t \le T \}$,
implies that
\begin{equs}
- J_\gamma (T,\rho) \le
\liminf_{v \rightarrow \infty}
\frac{1}{v}\log\Big(
{
\sup_{\|x_0^v\|_1 \le \gamma}} \,
\mathbb P_{x_0^v}\big[\sup_{t \in [0,T]} {\|X_t^v\|_1} > \rho \big]\Big) \,.
\label{e:lbd-sigma-rho}
\end{equs}
While proving \lref{t:exptight} we saw that the
\abbr{rhs} of \eref{e:lbd-sigma-rho} is, for some finite
$\kappa = \kappa(\gamma)$, with the constant $b$ of \aref{a:1}(a), any $T$ and $\rho \ge \rho(\ell)$,
at most
\begin{equs}
\limsup_{v \to \infty} v^{-1} \log \Big\{
\ell^{-v} [e^{\kappa v} + T e^{b v}] \Big\}
=  - \log \ell + \kappa \vee b \,.
\label{e:lbd-exptight}
\end{equs}
Combining \eref{e:lbd-sigma-rho} and \eref{e:lbd-exptight},
we establish \eref{e:quas-pot-inf}
upon taking $\rho \to \infty$ followed by $\ell \to \infty$.
\end{proof}

\section*{Acknowledgements}
We thank Ofer Zeitouni for pointing our attention to \cite{dupuis16}. Furthermore, AA and JPE would like to thank No\'{e} Cuneo and Neil Dobbs
for helpful discussions.

\bibliographystyle{JPE.bst}
\bibliography{bib.bib}

\end{document}